\documentclass[11pt]{amsart}
\usepackage{amsmath,amsfonts,amssymb,amscd}
\usepackage{latexsym}
\usepackage{multirow}
\usepackage[usenames]{color}
\usepackage[usenames,dvipsnames]{xcolor}
\usepackage{ifpdf}
\ifpdf
  \usepackage[pdftex]{graphicx}
  \usepackage{epstopdf}
\else
  \usepackage[dvips]{graphicx}
\fi
\usepackage{url}
\usepackage{pdflscape}
\usepackage[all]{xy}

\usepackage[top = 1in,bottom = 1in,left = 1in, right = 1in, foot = 1in,footnotesep=0.15in]{geometry}

\usepackage{setspace}

\newtheorem{thm}{Theorem}[section]

\theoremstyle{definition}

\newtheorem{conj}[thm]{Conjecture}

\usepackage[
hyperindex=true,pagebackref=true,bookmarks=true,
colorlinks=true,linkcolor=blue,citecolor=red]
{hyperref}
\def\~{{\rm --}}

\title [Exceptional Knot Homology]
{Exceptional Knot Homology}
\author[Ross Elliot]{Ross Elliot}
\author[Sergei Gukov]{Sergei Gukov}
\date{May 1, 2015}
\date{\today}

\address[R. Elliot]{
California Institute of Technology, Pasadena, California 91125, USA\\
relliot@caltech.edu\\}
\address[S. Gukov]{
California Institute of Technology, Pasadena, California 91125, USA\\
\newline
\indent \hspace{48pt}Simons Center for Geometry and Physics, Stony Brook, NY 11794, USA.}

 \def\bysame{{\bf --- }}
 \def\~{{\bf --}}
\newcommand{\comment}[1]{}
\renewcommand{\tilde}{\widetilde}
\renewcommand{\hat}{\widehat}

\begin{document}
\onehalfspacing
\begin{abstract}
The goal of this article is twofold.
First, we find a natural home for the double affine Hecke algebras (DAHA) in the physics of BPS states.
Second, we introduce new invariants of torus knots and links called \emph{hyperpolynomials}
that address the ``problem of negative coefficients'' often encountered in DAHA-based approaches
to homological invariants of torus knots and links.
Furthermore, from the physics of BPS states and the spectra of singularities associated with Landau-Ginzburg potentials,
we also describe a rich structure of differentials that act on homological knot invariants for exceptional groups
and uniquely determine the latter for torus knots. 
\end{abstract}

\maketitle

{\singlespacing
\tableofcontents }

\vfill\eject

\renewcommand{\natural}{\wr}

\setcounter{section}{-1}
\setcounter{equation}{0}
\onehalfspacing
\section{Introduction}

\comment{
We provide several examples of the hyperpolynomials
of type E6 for simplest torus knots and discuss
their symmetries and connections.
}

Categorification of quantum group invariants has been a very active area of research in the past several years.
By now, a number of methods have been developed that allow one to ``promote'' a polynomial invariant $P^{\frak g, V} (K;q)$
of a knot $K$ colored by a representation $V$ of $\mathcal{U}_q (\frak g)$
to a bi-graded homology theory $\mathcal{H}^{\frak g, V}_{i,j} (K)$,
whose Euler characteristic is $P^{\frak g, V} (K;q)$:
\begin{equation}
P^{\frak g, V} (K;q) \; = \; \sum_{i,j} (-1)^j q^i \text{dim} \mathcal{H}^{\frak g, V}_{i,j} (K) \,.
\label{PgR}
\end{equation}
In practice, it is often convenient to work with Poincar\'e polynomials of $\mathcal{H}^{\frak g, V}_{i,j} (K)$:
\begin{equation}
\mathcal{P}^{\frak g, V} (K;q,t) \; = \; \sum_{i,j} q^i t^j \text{dim} \mathcal{H}^{\frak g, V}_{i,j} (K) \,,
\label{PHgR}
\end{equation}
or, better yet, with the so-called superpolynomials $\mathcal{P} (K;a,q,t)$ that depend on three variables
and package homological invariants of arbitrary rank and fixed Cartan type.

While formal definitions of these homological knot invariants are available for many groups and representations \cite{KhR,Y,Web,Wu},
their calculation has been a daunting task. Besides the Khovanov-Rozansky homology \cite{KhR}, which corresponds to the fundamental representation of $\frak g = sl(N)$ and is reasonably computable, at present there exist only two approaches amenable to calculations for arbitrary groups and representations.
One approach \cite{DGR,GW,GS} is based on the formal
structure of knot homologies (and superpolynomials) that follows from the physical interpretation \cite{GSV,G,W2} of knot homologies.  Another approach \cite{ChJ,CJJ} proposed recently is based on DAHA (see also \cite{AS}).

Both of these approaches have advantages and disadvantages. The first approach allows one to compute homological invariants / superpolynomials of arbitrary knots, while the second approach
is limited to torus knots. On the other hand, the second approach can easily be implemented on a computer, whereas the first approach can only be done ``by hand'' for simple knots with up to 10 crossings or so.

More importantly for the present paper, neither approach is limited to classical groups or particular representations.
We use this feature to tackle one of the most difficult problems in this subject:
the study of homological invariants (and superpolynomials) associated with exceptional groups.
In fact, for a problem like this we will need to combine the power of both methods, because each one individually is not
sufficient for producing a polynomial with positive coefficients.

The simplest ``exceptional knot homology'' corresponds to the minuscule 27-dimensional representation of the simply-laced Lie algebra $\mathfrak{e}_6$. The representation of the principal $SL(2)$ on $\mathbf{27}$ is isomorphic to the representation
of the Lefschetz $SL(2)$ on the cohomology of the 16-dimensional flag variety $G/P$, with the Poincar\'e polynomial,

\footnotesize
\begin{equation*}
P(t)  =  1 + t + t^2 + t^3 + 2 t^4 + 2 t^5 + 2 t^6 + 2 t^7 + 3 t^8 
+ 2 t^9 + 2 t^{10} + 2 t^{11} + 2 t^{12} + t^{13} + t^{14} + t^{15} + t^{16}. 
\end{equation*}
\normalsize

The strategy of our approach will be the following.
First, we compute the DAHA-Jones polynomials of simple torus knots colored
by the 27-dimensional representation of $\mathfrak{e}_6$ using the approach of \cite{ChJ,CJJ}.
These will turn out to have both positive and negative coefficients.
To fix this problem and to construct analogues of superpolynomials with positive coefficients,
we will resort to the other method \cite{DGR,GS,GGS} based on a rich structure of the differentials.
Which differentials to expect and how they should act is controlled by deformations of
a certain singularity \cite{GW}, which will be yet another new result of this paper.
\medskip

The structure of this paper is as follows.  In Section \ref{HKBPS}, we will review the physical realization of knot homologies as spaces of BPS states in topological string theory.  In Section \ref{DAHAJ}, we define the DAHA-Jones polynomials and explain their relationship to torus knot polynomials and homologies.  

Section \ref{sec:hyper} contains our main proposal for $E_6$-hyperpolynomials, as well as three convincing examples.  At various intermediate stages in our calculations we shall need superpolynomials for root systems of Cartan types $A$ and $D$. The corresponding results are summarized in Appendix \ref{sec:AandD} and can be found in \cite{GW,GS,ChJ}.  Appendix \ref{fig} contains diagrams that depict our examples.

Finally, in Section \ref{sec:Adj} we classify the adjacencies (infinitesimal deformations) of the singularity $Z_{3,0}$ and compute the corresponding spectra.  As explained there and in \cite{GW}, deformations of this singularity control which differentials we are to expect from $(\mathfrak{e}_6,\mathbf{27})$ knot homologies.  The results of this analysis are contained in Appendix \ref{app:Adj}.

\subsection{Acknowledgements}
Our special thanks go to Ivan Cherednik, who provided the formulas for DAHA-Jones polynomials and participated in the development of many ideas contained herein.  Without his contributions, this work would not be possible.

We would also like to thank J.~Adams, M.~Aschbacher, D. Bar-Natan, P. Cvitanovi\'{c}, W.A. de Graaf, A.~Gabrielov, and S. Morrison
for helpful discussions.  The work of S.G. is funded in part by the DOE Grant DE-SC0011632 and the Walter Burke Institute for Theoretical Physics.  The work of R.E. is partially supported by a Troesh Family Graduate Fellowship 2014-15.
\medskip

\setcounter{equation}{0}

\section{Knot homologies and refined BPS states}\label{HKBPS}

\subsection{Large $N$ duality and BPS states}

Following \cite{W1}, recall that the Chern-Simons TQFT on a 3-manifold $M$ with gauge group $G$ at \emph{level} $k\in\mathbb{Z}$ is described by the action functional:
\begin{equation}
S(A) =  \frac{k}{4\pi}\displaystyle\int_M \text{tr} (A\wedge dA + \frac{2}{3}A\wedge A\wedge A),
\end{equation}
where $A$ is the ($\mathfrak{g}$-valued) connection one-form of a principal $G$-bundle on $M$.  The partition function of this theory is given by the path integral,
\begin{equation}
Z(M) = \displaystyle\int_\mathcal{A}[\mathcal{D}A]e^{iS(A)},
\end{equation}
over the configuration space $\mathcal{A}$ of principal $G$-connections on $M$.  Owing to the topological nature of Chern-Simons theory, $Z(M)$ is, \emph{a fortiori}, a topological invariant of $M$.

Now consider the open string theory described by the topological A-model on the cotangent bundle $T^*M$ with $N$ D-branes wrapping the Lagrangian $M\subset T^*M$, and coupling constant,
\begin{equation}
g_s = \frac{2\pi i}{k+N}.
\end{equation}
When $G = SU(N)$, it was shown in \cite{W3} that the $\frac{1}{N}$ expansion of the Chern-Simons free energy $F(M) = \log Z(M)$ is naturally identified with with the contribution to free energy of the degenerate instantons in this topological string setup.

Instantons there are generally described by holomorphic maps of Riemann surfaces with Lagrangian boundary  conditions:
\begin{equation}
(\Sigma,\partial\Sigma)\hookrightarrow(T^*M,M).
\end{equation}
However, an easy consequence of Witten's ``vanishing theorem" is that the \emph{only} such maps are the degenerate (constant) ones.  Therefore, one identifies
\begin{equation}
Z_{\text{CS}}(M) = Z^{\text{open}}_{\text{string}}(T^*M),
\end{equation}
the partition functions for Chern-Simons gauge theory on $M$ and the open topological string theory on $T^*M$.

In the special case of $M = {\bf S}^3$, it was conjectured \cite{GV} that at large $N$, this open string setup undergoes a geometric transition which produces a (physically equivalent) closed string theory.  This \emph{conifold transition} shrinks the 3-cycle of the \emph{deformed conifold} $T^*{\bf S}^3$ to a point and resolves the resulting conical singularity with a small blow-up.  The resulting space $X$ is the \emph{resolved conifold}, i.e. the total space of the $\mathcal{O}(-1)\oplus\mathcal{O}(-1)$ bundle over $\mathbb{C}{\bf P}^1$.

Observe that the conifold transition eliminates the $N$ branes wrapping ${\bf S}^3$, producing a closed string theory on $X$.  In the worldsheet description of this theory, based on the genus $g$ topological sigma model coupled to 2-dimensional gravity, the free energy is
\begin{equation}
F_g(t) = \sum_{Q\in H_2(X)}N_{g,Q}e^{-tQ},
\end{equation}
where the parameter $t$ is the K$\ddot{\text{a}}$hler modulus for the Calabi-Yau space $X$:
\begin{equation}
t = \frac{2\pi iN}{k+N} = \text{vol}(\mathbb{C}{\bf P}^1),
\label{tVolN}
\end{equation}
and $N_{g,Q}$ is the Gromov-Witten invariant ``counting" holomorphic maps of genus $g$ representing the integral 2-homology class $Q$.

The numbers $N_{g,Q}$ are rational, in general.  However, as shown in \cite{GV}, this model also admits a target space description in which the all-genus free energy is naturally described in terms of integer invariants $n_Q^s\in\mathbb{Z}$:
\begin{equation}
F(g_s,t) = \displaystyle\sum_{g=0}^\infty g_s^{2g-2}F_g(t) = \displaystyle\sum_{\substack{Q\in H_2(X),\\ s\geq 0}}n_Q^s\left[\sum_{m\geq 1}\frac{1}{m}\left(2\sin \frac{mg_s}{2} \right)^{2s-2}e^{-mtQ}\right],
\end{equation}
which encode degeneracies of the so-called \emph{BPS states}.

In a general supersymmetric quantum theory, a BPS state is one whose mass is equal to the central charge of the supersymmetry algebra.  In the case at hand, a state is a D2-brane wrapping $\mathbb{C}{\bf P}^1$, and the BPS condition means that it is supported on a a calibrated 2-submanifold of the Calabi-Yau $X$ (i.e. on a holomorphic curve in $X$).

Thus, a minimally embedded surface representing $Q\in H_2(X;\mathbb{Z})$ gives rise to a component of the Hilbert space $\mathcal{H}_{BPS}$, i.e. a projective unitary representation of the spatial rotation group,
\begin{equation}
SO(4) \sim SU(2)_L \times SU(2)_R,
\end{equation}
of $\mathbb{R}^4$ obtained upon compactification from M-theory. This representation can be specified by two half-integer charges $j_L,j_R\in\frac{1}{2}\mathbb{Z}_{\geq 0}$, which are the weights of the respective $SU(2)$ representations.

One might be tempted to introduce integers $n_Q^{(j_L,j_R)}$ counting these states.  However as one deforms the theory, BPS states can combine into non-BPS states, so these numbers are not invariant.  On the other hand, the index,
\begin{equation}
n_Q^{j_L} := \displaystyle\sum_{j_R} (-1)^{2j_R}(2j_R+1)n_Q^{(j_L,j_R)},
\end{equation}
is well-defined on the moduli of $X$.  The integers $n_Q^s$ are then related by a change of basis for the representation ring of $SU(2)$.

\subsection{Knot invariants and topological strings}

For a knot $K \subset M$ and a representation $V$ of $\mathfrak{g}$, one can consider the holonomy of $A$ along $K$ traced in $V$, yielding the gauge-invariant Wilson loop operator:
\begin{equation}
W_V^K(A) = \text{tr}_V\left[\mathcal{P}\text{exp}\displaystyle\oint_K A\right].
\end{equation}
Expanding the correlation function of a Wilson loop in $q := e^\frac{2\pi i}{k+h^\vee}$ produces an integer Laurent polynomial:
\begin{equation}
P^{\frak g, V} (M,K;q) := \left\langle W_V^K\right\rangle_M = \frac{1}{Z(M)}\displaystyle\int_\mathcal{A}[\mathcal{D}A]e^{iS(A)}W_V^K(A),
\end{equation}
which is naturally an isotopy invariant of $K\subset M$.  In what follows, we will exclusively consider $K\subset {\bf S}^3$ and suppress $M$.  Then $P^{\frak g, V} (K;q)$ are the \emph{quantum knot invariants} discussed in the introduction and whose categorifications \eqref{PgR} we will discuss below.

As explained in \cite{OV}, Wilson loops can be incorporated in the open string on the deformed conifold by introducing $L_K \subset T^* {\bf S}^3$, the \emph{conormal bundle} to $K\subset {\bf S}^3$.  In particular, $L_K$ is a Lagrangian submanifold of $T^* {\bf S}^3$, which is topologically ${\bf S}^1 \times \mathbb{R}^2$ and with $L_K \cap {\bf S}^3 = K$.  Wrapping $M$ ``probe" branes on $L_K$ produces a theory with three kinds of strings:
\begin{enumerate}
\item both ends on ${\bf S}^3$ $\leadsto$ $SU(N)$ Chern-Simons theory on ${\bf S}^3$,
\item both ends on $L_K$ $\leadsto$ $SU(M)$ Chern-Simons theory on $L_K$,
\item one end on each ${\bf S}^3$ and $L_K$ $\leadsto$ complex $SU(N)\otimes SU(M)$ scalar field on $K$.
\end{enumerate}
Let $U$, $V$ be the holonomies around $K$ of gauge fields $A$, $A'$ in (1),(2) respectively.  Then the last kind of string (3) contributes to the overall action by
\begin{equation}
S(U,V) := \displaystyle\sum_{n=1}^\infty\frac{1}{n}\text{tr}U^n\text{tr}V^{-n} = \log\left[\displaystyle\sum_R\text{tr}_RU\text{tr}_RV^{-1}\right].
\end{equation}
In turn, the effective action for the theory on ${\bf S}^3$ is
\begin{equation}
S(A; K) := S_{CS}(A;{\bf S}^3) + S(U,V),
\end{equation}
and integrating $A$ out of the overall theory involves evaluating
\begin{equation}
\langle S(U,V)\rangle_{{\bf S}^3} = \frac{1}{Z({\bf S}^3)}\displaystyle\int_\mathcal{A}[\mathcal{D}A]e^{iS(A; K)} = \sum_\lambda\langle W_\lambda^K\rangle(\text{tr}_\lambda V^{-1})
\end{equation}
for fixed $V$, which produces a generating functional for all Wilson loops associated to $K\subset {\bf S}^3$ (i.e. for all Young diagrams $\lambda$).

If one follows the Lagrangian $L_K \subset T^*{\bf S}^3$ through the conifold transition, the result is another Lagrangian $L_K'\subset X$, where the $M$ branes will still reside.   In the resulting open string theory, the worldsheet perspective again ``counts", in an appropriate sense, holomorphic maps of Riemann surfaces with Lagrangian boundary conditions:
\begin{equation}
(\Sigma,\partial\Sigma)\hookrightarrow(X,L_K'),
\end{equation}
described by the open Gromov-Witten theory.

From the target space perspective, states correspond to configurations in which D2-branes wrap relative cycles $Q\in H_2(X,L_K';\mathbb{Z})$ and end on D4-branes which wrap $L_K'$.  BPS states are then minimally-embedded surfaces $\Sigma \subset X$ with boundaries $\partial\Sigma\subset L_K'$.

In \cite{OV}, the authors also showed that the generating functional for Wilson loops has an interpretation in terms of BPS degeneracies:
\begin{equation}
\langle S(U,V)\rangle_{{\bf S}^3} =  i\displaystyle\sum_{R,Q,s}N_{R,Q,s}\left[\sum_{m\geq 1}\frac{e^{m(-tQ+isg_s)}}{2m\sin\left(\frac{mg_s}{2}\right)}\text{tr}_RV^m\right],
\end{equation}
where $N_{R,Q,s}\in\mathbb{Z}$ are certain modifications of $n_Q^s$.  One can then express the quantum invariant $P^{\frak{sl}_N, R} (K;q)$ directly in these terms.  For example, if $R=\square$ we have
\begin{equation}
P_N (K;q) = \frac{1}{q-q^{-1}}\displaystyle\sum_{Q,s}N_{\square,Q,s}q^{NQ+s},
\end{equation}
directly relating quantum knot invariants to the enumerative geometry of $X$.

\subsection{Knot homologies and refined BPS states}

In light of the mathematical development of homology theories categorifying quantum knot invariants, one might ask whether they also admit physical descriptions in the contexts outlined above.  This program was initiated in \cite{GSV}, where the authors refined the BPS degeneracies:
\begin{equation}
N_{\square,Q,s}(K) = \displaystyle\sum_r(-1)^rD_{Q,s,r}(K),
\end{equation}
introducing \emph{non-negative} integers $D_{Q,s,r}\in\mathbb{Z}_{\geq 0}$, which also reflect the charge $r$ of $U(1)_R\in SU(2)_R$.
Given that the Calabi-Yau $X$ is rigid, these numbers are invariant under complex structure
deformations.\footnote{Furthermore, as we mentioned earlier in \eqref{tVolN}, the K\"ahler modulus of $X$ is related to the rank
of the underlying root system via $q^N = e^t = \exp ( \text{vol}(\mathbb{C}{\bf P}^1))$,
so that changes in the BPS spectrum as one varies the K\"ahler parameter $t$ (a.k.a. the `stability parameter')
reflect changes of homological knot invariants at different values of $N$. See \cite{GS} for details.}

This led to a conjecture relating the knot homology which categorifies $P_N(K;q)$ to refined BPS degeneracies:
\begin{equation}
(q-q^{-1})KhR_N(K; q,t) = \displaystyle\sum_{Q,s,r}D_{Q,s,r}(K)q^{NQ+s}t^r,
\end{equation}
for sufficiently large $N$, where $KhR_N(K; q,t)$ is the Poincar\'{e} polynomial for the Khovanov-Rozansky homology.

More generally, one might view the charges $Q,s,r$ as gradings on the Hilbert space $\mathcal{H}_{\text{BPS}}(K)$ and conjecture an isomorphism of graded vector spaces:
\begin{equation}
\displaystyle\bigoplus_{i,j}H_{i,j}(K) = H_{\text{knot}}(K) \cong \mathcal{H}_{\text{BPS}}(K) = \displaystyle\bigoplus_{Q,s,r} \mathcal{H}_{Q,s,r}(K),
\label{HspacesPh}
\end{equation}
with $\dim\mathcal{H}_{Q,s,r}(K) =  D_{Q,s,r}(K)$.   This new perspective has revealed hidden structures of knot homologies that are manifest in the context of BPS states.  In particular, $H_{\text{knot}}(K)$ should:
\begin{itemize}
\item stabilize in dimension for sufficiently large $N$
\item be triply-graded, with the additional grading (corresponding to $Q$) encoding $N$-dependence of the homology theory
\item include the structure of differentials (c.f. Section \ref{approach}) corresponding to wall-crossing behavior of $\mathcal{H}_{\text{BPS}}(K)$
\end{itemize}
and, in fact, all of these structures were realized in \cite{DGR}, where the authors proposed a triply-graded homology theory categorifying the HOMFLY polynomial.  Furthermore, they were able to construct explicit Poincar\'{e} polynomials for this homology theory (``superpolynomials") based on a rigid structure of differentials, which was later formalized in \cite{R}.  Similar constructions for other choices of $(\mathfrak{g},V)$ were proposed in \cite{GW,GS,GGS}.

\subsection{M-theory descriptions}

M-theory on an eleven-dimensional space-time incorporates the various (equivalent) versions of string/gauge theory and the dualities between them.  The individual theories can then be recovered by integrating out the dependence of M-theory on some portion of the background geometry.

Naturally, this framework can offer several equivalent but nontrivially different points of view on the same object.  In the case of knot homologies, we are looking for new descriptions of
\begin{equation}
H_{\text{knot}}(K) \cong \mathcal{H}_{\text{BPS}}(K),
\end{equation}
so we promote the topological string setups described above.

In particular, the five-brane configuration relevant to the physical description of the $(\mathfrak{sl}_N, \lambda)$ knot homologies on the deformed conifold is:
\begin{eqnarray}
\text{space-time} & : & \qquad \mathbb{R} \times T^* {\bf S}^3 \times M_4 \nonumber \\
N~\text{M5-branes} & : & \qquad \mathbb{R} \; \times \; {\bf S}^3 \; \times \; D \\
|\lambda|~\text{M5-branes} & : &  \qquad \mathbb{R} \; \times \; L_K \; \times \; D \nonumber
\label{M51}
\end{eqnarray}
and the equivalent (large-$N$ dual) configuration on the resolved conifold is:
\begin{eqnarray}
\text{space-time} & : & \qquad \mathbb{R} \; \times \; X \; \times \; M_4 \label{M52} \\
|\lambda|~\text{M5-branes} & : &  \qquad \mathbb{R} \; \times \; L_K' \; \times \; D \nonumber
\end{eqnarray}
where states correspond to configurations in which M2-branes wrap relative cycles $Q\in H_2(X,L_K';\mathbb{Z})$, fill $\mathbb{R}$, and end on the M5-branes.

The precise form of the 4-manifold $M_4$ and the surface $D \subset M_4$ is not important (in most applications $D \cong \mathbb{R}^2$ and $M_4 \cong \mathbb{R}^4$),
as long as they enjoy a $U(1)_F \times U(1)_P$ symmetry action, corresponding to the charges that comprise the $(s,r)$-gradings. The first (resp. second) factor is a rotation symmetry of the normal (resp. tangent) bundle of $D \subset M_4$.  Following \cite{W2}, let us denote the corresponding quantum numbers by $F$ and $P$.
These quantum numbers were denoted, respectively, by $2 S_1$ and $2(S_1 - S_2)$ in \cite{AS} and by $2 j_3$ and $n$ in \cite{GS}.

This description of $\mathcal{H}_{\text{BPS}}(K)$ in the M-theory framework led to a number of developments
which shed light on various aspects of knot homologies and yield powerful computational techniques.
Some examples include:
\begin{itemize}
\item \cite{W2} formulates the relevant space of BPS states within (\ref{M52})
\item \cite{AS} refines torus knot invariants directly within Chern-Simons theory based on its relationship with (\ref{M52}) discovered in \cite{W3}
\item \cite{DGH} takes the perspective of $M_4$ on which the BPS invariants are expressed via equivariant instanton counting
\end{itemize}
\medskip

\setcounter{equation}{0}
\section{DAHA-Jones polynomials}\label{DAHAJ}

Given the $(r,s)$-torus knot, a root system $R$, and a weight $b$, the corresponding \emph{DAHA-Jones polynomial} is defined by the simple formula:
$$JD^R_{r,s}(b; q,t) :=\{\hat{\gamma}_{r,s}(P_b)/P_b(q^{-\rho_k})\}_{ev}$$
We will briefly explain the meaning of this expression and then describe its properties and relations to torus knot polynomials and homologies.

Good general references for the material in this section are \cite{C101,Ha,Hu,Ki,M1,M5} as well as the original papers \cite{C15,ChJ,CJJ,M2,M3}.  Our conventions for root systems will be from \cite{Bo}.

\subsection{Affine Hecke algebras}

\subsubsection{\sf Hecke algebras}

Let $R$ be a (crystallographic) root system of rank $n$ with respect to the Euclidean inner product $(-,-)$ on $\mathbb{R}^n$, and let $\Delta = \{\alpha_1,\ldots,\alpha_n\}$ be any set of simple roots.  The \emph{Weyl group} $W$ for $R$ is generated by the simple reflections:
\begin{equation}
s_i: \beta\mapsto \beta - \frac{2(\beta,\alpha_i)}{(\alpha_i,\alpha_i)}\alpha_i\hspace{5pt}\text{for}\hspace{5pt} 1\leq i\leq n, \hspace{5pt}\beta \in R,
\end{equation}
subject to the Coxeter relations $(s_is_j)^{m_{ij}} = 1$.  The numbers $m_{ij}$ are 2,3,4,6 when the corresponding nodes in the Dynkin diagram for $R$ are joined by 0,1,2,3 edges, respectively.

The (nonaffine) \emph{Hecke algebra} $H$ for $R$ is generated over $\mathbb{C}(t_1^{\frac{1}{2}},\ldots,t_n^{\frac{1}{2}})$ by elements $\{T_1,\ldots,T_n\}$, subject to relations:
\begin{equation}
(T_i - t_i^{\frac{1}{2}})(T_i+t_i^{-\frac{1}{2}})\text{ for }1\leq i\leq n,
\end{equation}
\begin{equation}
T_iT_jT_i\ldots = T_jT_iT_j\ldots\text{ with }m_{ij}\text{ terms on each side},
\end{equation}
where the number of distinct $t_i$ is equal to the number the orbits of $W$ acting on $R$, so at most 2 in the nonaffine case.  That is, we normalize the form by $(\alpha,\alpha)=2$ for short roots $\alpha\in R$ and set $\nu_\beta :=  \frac{(\beta,\beta)}{2}$ for $\beta\in R$.  Then $t_i :=  t_{\nu_{\alpha_i}}$ for each simple root $\alpha_i\in\Delta$.

\subsubsection{\sf Twisted affine root systems}

Before defining an affine root system, we recall the identification $\mathbb{R}^{n+1} \cong \text{Aff}(\mathbb{R}^n)$.  That is, we interpret a vector $[\vec{u},c]\in\mathbb{R}^n\times\mathbb{R}$ as an affine linear function on $\mathbb{R}^n$:
\begin{equation}
[\vec{u},c]: \vec{v}\mapsto (\vec{u},\vec{v}) - c,
\end{equation}
whose zero set $[\vec{u},c]^{-1}(0)$ is an affine hyperplane in $\mathbb{R}^n$, $H_{[\vec{u},c]} :=  \{\vec{v}\in\mathbb{R}^n : (\vec{u},\vec{v})=c\}$.  Observe that $H_{[\vec{u},c]} = H_{[\vec{u},0]} + \frac{c}{2}\vec{u}^\vee$, where $\vec{u}^\vee:=  \frac{\vec{u}}{\nu_{\vec{u}}}$.  

The reflection of $\mathbb{R}^n$ through $H_{[\vec{u},c]}$ is
\begin{equation}
s_{[\vec{u},c]} : \vec{v} \mapsto \vec{v} - \left[(\vec{u},\vec{v}) - c\right]\vec{u}^\vee,
\end{equation}
which fixes $H_{[\vec{u},c]}$ and maps 0 to $c\vec{u}^\vee$.  We can extend the domain of affine reflections to act on $\text{Aff}(\mathbb{R}^n) \cong \mathbb{R}^n\times\mathbb{R}$ by
\begin{equation}
s_{[\vec{u},c]}([\vec{v},k]) :=  [\vec{v},k]\circ s_{[\vec{u},c]} = [\vec{v},k] - (\vec{v},\vec{u}^\vee)[\vec{u},c].
\end{equation}
Alternatively, we could describe $s_{[\vec{u},c]}$ as a reflection in $H_{[\vec{u},0]}$ with a subsequent translation by $c\vec{u}^\vee$, where ``translations" are
\begin{equation}
s_{[\pm\vec{u},c]}s_{[\vec{u},0]} = s_{[\vec{u},0]}s_{[\mp\vec{u},c]}\text{ : }
\begin{cases}
\vec{v}\mapsto\vec{v}\pm c\vec{u}^\vee,\\
[\vec{v},k]\mapsto[\vec{v},k\pm(\vec{v},\vec{u}^\vee)c],
\end{cases}
\end{equation}
and we will often confuse $c\vec{u}^\vee\in\mathbb{R}^n$ with this action below.

Define the \emph{(twisted) affine root system} $R\subset\widetilde{R}$ by:
\begin{equation}
\widetilde{R} = \{[\alpha,k\nu_\alpha] : \alpha\in R, k\in\mathbb{Z}\},
\end{equation}
with $R = \{[\alpha,0]\}$.  The simple roots for $\widetilde{R}$ are $\widetilde{\Delta}:=  \{\alpha_0 = [-\vartheta,1]\}\cup\Delta$, where $\vartheta\in R$ is the highest \emph{short} root with respect to $\Delta$.

\subsubsection{\sf Affine Weyl groups}

The \emph{affine Weyl group} $\widetilde{W}$ is generated by $s_i: = s_{\alpha_i}$, $0\leq i\leq n$ subject to relations $s_i^2 = 1$ and
\begin{equation}
s_is_js_i\ldots = s_js_is_j\ldots\text{ with }m_{ij}\text{ terms on each side},
\end{equation}
where $m_{ij}$ correspond, as above, to the affine Dynkin diagram.  

We saw that $s_{[\alpha,k\nu_{\alpha}]}$ admits a description as a reflection $s_\alpha\in W$ composed with a translation by $k\nu_{\alpha}\alpha^\vee = k\alpha\in Q$, where $Q$ is the root lattice for $R$, i.e., the $\mathbb{Z}$-span of $\Delta$.  Therefore, one easily concludes that
\begin{equation}
\widetilde{W} = W \ltimes Q,
\end{equation}
 where $Q$ acts by ``translations" as described above.

If we enlarge the group $Q$ to include translations by the weight lattice,
\begin{equation}
Q\subset P:=  \displaystyle\bigoplus_{i=1}^n\mathbb{Z}\omega_i,
\end{equation}
where $\{\omega_i\}$ are fundamental weights, we obtain the \emph{extended affine Weyl group},
\begin{equation}
\widehat{W}:=  W \ltimes P = \widetilde{W}\ltimes\Pi,
\end{equation}
where $\Pi:= P/Q$ in the semidirect product decomposition relative to $\widetilde{W}$.  

To describe the subgroup $\Pi\triangleleft \widehat{W}$ more explicitly, we can introduce a length function $l$ on $\widehat{W}$:
\begin{equation}
l(\hat{w}):=  \left|\widetilde{R}_+ \cap \hat{w}(-\widetilde{R}_+)\right|,
\end{equation}
where $\widehat{R}_+$ is the set of positive roots with respect to $\widetilde{\Delta}$.  Then $\Pi = \{\hat{w}\in\widehat{W} : l(\hat{w})=0\}$. Geometrically, these these are the elements of $\widehat{W}$ which permute $\widetilde{\Delta}$, and we can label an element $\pi_r\in\Pi$ by its action $\pi_r(\alpha_0) = \alpha_r$.

Alternatively, define the set of indices of minuscule weights:
\begin{equation}
O' :=  \{r : 0\leq (\omega_r, \alpha^\vee)\leq 1\text{, }\text{, for all }\alpha\in R_+\}\subset\{1,\ldots,n\}.
\end{equation}
Then $O = \{0\}\cup O'$ is a system of representatives for $P/Q$ in the sense that every $b\in P$ can be written uniquely as $b = \omega_r + \alpha$ for some $r\in O$, $\alpha\in Q$, where $\omega_0 = 0$.  For $r\in O$ let $u_r\in W$ be the shortest element such that $u_r(\omega_r)\in -P_+$.  We can define
\begin{equation}
\Pi = \{\pi_r : \omega_r = \pi_ru_r\text{, }r\in O\},
\end{equation}
and observe that $\pi_0 = \text{id}$.  

The affine Weyl group $\widetilde{W}$ (or, to be more precise, its group algebra) has a simple physical interpretation \cite{GW2}
as the algebra of line operators in four-dimensional gauge theory on $M_4 \cong {\bf S}^1 \times \mathbb{R}^3$
in the presence of ramification along $D \cong {\bf S}^1 \times \mathbb{R}$. (In physics, ramification is often called a \emph{surface operator}.)

\subsubsection{\sf Affine Hecke algebras}
The \emph{affine Hecke algebra} $\mathcal{H}$ for $R\subset\widetilde{R}$ is generated over $\mathbb{C}(t_0^{\frac{1}{2}},t_1^{\frac{1}{2}},\ldots,t_n^{\frac{1}{2}})$.  It admits two equivalent descriptions, each emphasizing one of the two equivalent descriptions of the extended affine Weyl group $\widehat{W}$:  
\begin{itemize}
\item For $\widehat{W} = \widetilde{W}\ltimes\Pi$, $\mathcal{H}$ is generated by elements $\{T_0,T_1,\ldots,T_n\}$ and $\pi_r\in\Pi$, subject to relations:
\begin{enumerate}
\item $(T_i - t_i^{\frac{1}{2}})(T_i+t_i^{-\frac{1}{2}})$ for $0\leq i\leq n$,
\item $T_iT_jT_i\ldots = T_jT_iT_j\ldots$ with $m_{ij}$ terms on each side,
\item $\pi_rT_i\pi_r^{-1} = T_j$ if $\pi_r(\alpha_i) = \alpha_j$.
\end{enumerate}

\item For $\widehat{W} = W\ltimes P$, $\mathcal{H}$ is generated by $\{T_1,\ldots,T_n\}$ and $\{Y_b : b\in P\}$, subject to relations:
\begin{enumerate}
\item $(T_i - t_i^{\frac{1}{2}})(T_i+t_i^{-\frac{1}{2}})$ for $1\leq i\leq n$,
\item $T_iT_jT_i\ldots = T_jT_iT_j\ldots$ with $m_{ij}$ terms on each side,
\item $Y_{b+c} = Y_bY_c\text{ for }b,c\in P$,
\item $T_iY_b = Y_bY_{\alpha_i}^{-1}T_i^{-1}\text{ if }(b,\alpha^\vee_i) = 1\text{ for }0\leq i\leq n$,
\item $T_iY_b = Y_bT_i\text{ if }(b,\alpha^\vee_i) = 0\text{ for }0\leq i\leq n$.
\end{enumerate}
\end{itemize}

To translate from the first to the second description, one can define pairwise-commuting elements:
\begin{equation}\label{Yops}
Y_b:= \displaystyle\prod_{i=1}^nY_i^{l_i}\text{ for }b = \displaystyle\sum_{i=1}^nl_i\omega_i\in P,
\end{equation}
where $Y_i:=  T_{\omega_i}$ for $\omega_i\in\widehat{W}$.  That is, if $l=l(\tilde{w})$ so that $\tilde{w} = s_{i_l}\cdots s_{i_1}\in\widetilde{W}$ is a reduced decomposition, then $T_{\pi_r\tilde{w}}:= \pi_rT_{i_l}\cdots T_{i_1}$.  For example, $Y_\vartheta = T_0T_{s_\vartheta}$.

Much like the affine Weyl group, the affine Hecke algebra $\mathcal{H}$ can also be interpreted as
the algebra of line operators in 4d gauge theory on $M_4$ with a ramification (surface operator) along $D \subset M_4$.
The only difference is that now one has to introduce the so-called $\Omega$-background in the normal bundle of $D$.
(See \cite{G2} for a review.)

\subsection{DAHA and Macdonald polynomials}
\subsubsection{\sf Double affine Hecke algebras}

Let $m$ be the least natural number satisfying $(P,P) \subset  \frac{1}{m} \mathbb{Z}$.  Suppose that $\tilde{b} = [b,j]$ with $b = \displaystyle\sum_{i=1}^nl_i\omega_i \in P$ and $j\in\frac{1}{m}\mathbb{Z}$.  Then for $\{X_1,\ldots,X_n : [X_i,X_j]=0\}$ we define elements:
\begin{equation}
X_{\tilde{b}} :=  \displaystyle\prod_{i=1}^n X_i^{l_i}q^j,
\end{equation}
and an action of $\hat{w}\in\widehat{W}$ by $\hat{w}(X_{\tilde{b}}):= X_{\hat{w}(\tilde{b})}$.  Observe that $X_0:= X_{\alpha_0} = qX_\vartheta^{-1}$.

The \emph{double affine Hecke algebra} (``DAHA") $\mathcal{H\!\!\!H}$ for $R\subset\widetilde{R}$ is generated over $\mathbb{Z}_{q,t}:= \mathbb{Z}[q^{\frac{1}{m}},t_\nu^{\frac{1}{2}}]$
by elements $\{T_i, X_b, \pi_r : 0\leq i\leq n\text{, }b\in P\text{, } r\in O\}$ subject to relations:
\begin{enumerate}
\item $(T_i - t_i^{\frac{1}{2}})(T_i+t_i^{-\frac{1}{2}})$ for $0\leq i\leq n$,
\item $T_iT_jT_i\ldots = T_jT_iT_j\ldots$ with $m_{ij}$ terms on each side,
\item $\pi_rT_i\pi_r^{-1} = T_j$ if $\pi_r(\alpha_i) = \alpha_j$,
\item $T_iX_b = X_bX_{\alpha_i}^{-1}T_i^{-1}\text{ if }(b,\alpha^\vee_i) = 1\text{ for }0\leq i\leq n$,
\item $T_iX_b = X_bT_i\text{ if }(b,\alpha^\vee_i) = 0\text{ for }0\leq i\leq n$,
\item $\pi_rX_b\pi_r^{-1} = X_{\pi_r(b)} = X_{u_r^{-1}(b)}q^{(\omega_{\iota(r)},b)}\text{ for }r\in O'$,
\end{enumerate}
where in (6) we have used the involution $\iota:O'\rightarrow O'$ defined by $\pi_r^{-1} = \pi_{\iota(i)}$.

Observe that $\mathcal{H\!\!\!H}$ contains two subalgebras isomorphic to the affine Hecke algebra $\mathcal{H}$ for $R\subset\widetilde{R}$:
\begin{gather}
\mathcal{H}_1 :=  \langle\pi_r,T_0,\ldots,T_n\rangle\subset\mathcal{H\!\!\!H},\\
\mathcal{H}_2 :=  \langle T_1,\ldots,T_n,X_b\rangle\subset\mathcal{H\!\!\!H}.
\end{gather}
One can make $\mathcal{H}_1$ look more like $\mathcal{H}_2$ by defining pairwise-commuting elements $Y_b$ as in \eqref{Yops}.  Then we have that
\begin{equation}
\mathcal{H}_1 = \langle T_1,\ldots,T_n, Y_b\rangle.
\end{equation}
In fact, $\mathcal{H\!\!\!H}$ is also generated by elements $\{X_a,T_w,Y_b : a,b\in P\text{, }w\in W\}$.  While relations between these generators are more complicated, this presentation has some nice properties that will be useful in our definitions of Macdonald and DAHA-Jones polynomials below.  In particular, we have the PBW theorem for DAHA.

\begin{thm}\label{PBW} ($\text{PBW Theorem}$)
Any $h\in\mathcal{H\!\!\!H}$ can be written uniquely in the form
\begin{equation}
h = \displaystyle\sum_{a,w,b}c_{a,w,b}X_a T_w Y_b,
\end{equation}
for $c_{a,w,b}\in\mathbb{Z}_{q,t}$.  The similar statement holds for each ordering of $\{X_a,T_w,Y_b\}$.
\end{thm}

Much like the affine Weyl group and the affine Hecke algebra, the double affine Hecke algebra $\mathcal{H\!\!\!H}$
can be interpreted as the algebra of line operators in the presence of ramification (surface operator) \cite{G2}.

\subsubsection{\sf Polynomial representation}
To define the Macdonald polynomials using DAHA, we need the \emph{polynomial representation} 
\begin{equation}
\varrho: \mathcal{H\!\!\!H}\rightarrow\mathcal{V},
\end{equation}
where $\mathcal{V} :=  \text{End}(\mathbb{Z}_{q,t}[X])$. In generators $\{X_b,\pi_r,T_i\}$ its action is given by
\begin{equation}
 \varrho:
  \begin{cases}
   X_b\cdot g = X_bg\\
   \pi_r\cdot g = \pi_rg\pi_r^{-1}\text{, where, e.g., }\pi_r\cdot X_b = X_{\pi_r(b)}\\
   T_i\cdot g = \widehat{T}_ig
  \end{cases},
\end{equation}
for $g\in\mathbb{Z}_{q,t}[X]$.  The action of $T_i$ is by the \emph{Demazure-Lusztig} operators:
\begin{equation}\label{DL}
\widehat{T}_i :=  t_i^{\frac{1}{2}}s_i + (t_i^{\frac{1}{2}}-t_i^{\frac{-1}{2}})\frac{s_i-1}{X_{\alpha_i}-1},
\end{equation}
where, again, $s_iX_b = X_{s_i(b)}$.  Observe that if $g\in\mathbb{Z}_{q,t}[X]^W$ is \emph{any} symmetric polynomial, then $\widehat{T}_ig = t_i^{\frac{1}{2}}g$.  Remarkably, $\varrho$ is a faithful representation.

\subsubsection{\sf Symmetric Macdonald polynomials}

The symmetric Macdonald polynomials $P_b \in \mathbb{Z}_{q,t}[X]$ for $b\in P_+$ were introduced in \cite{M2,M3}.  They form a basis for the symmetric ($W$-invariant) polynomials $\mathbb{Z}_{q,t}[X]^W$.  DAHA provides a uniform construction of $P_b$ for any root system as the simultaneous eigenfunctions for a commuting family of $W$-invariant operators $L_f$ for $f\in\mathbb{Z}_{q,t}[Y]^W = Z(\mathcal{H}_1)$, see \cite{C15}.

Now for $f\in\mathbb{Z}_{q,t}[Y]^W\subset\mathcal{H\!\!\!H}$, we can use the polynomial representation to write an operator $L_f:= \varrho(f)$ on $\mathbb{Z}_{q,t}[X]$.
The \emph{symmetric Macdonald polynomials} are uniquely defined by
\begin{equation}
L_f(P_b) = f(q^{\rho_k+b})P_b,
\end{equation}
as simultaneous eigenfunctions of the pairwise-commuting $W$-invariant operators $L_f$ for all $f\in\mathbb{Z}_{q,t}[Y]^W$.
 In fact, $P_b\in\mathbb{Q}(q,t_\nu)[X]^W$.

In expressing $P_b$ as an eigenfunction, we used the notation \begin{equation}
\rho_k :=  \frac{1}{2}\displaystyle\sum_{\alpha\in R_+} k_\alpha\alpha = k_{\text{sht}}\rho_{\text{sht}} + k_{\text{lng}}\rho_{\text{lng}}\text{, }\hspace{5pt}\text{where, e.g., }\hspace{5pt}\rho_{\text{sht(lng)}} :=  \frac{1}{2}\sum_{\substack{\alpha \text{ short}\\\text{(long)}}} k_\alpha\alpha,
\end{equation}
for the Weyl vector weighted by a function $k_\alpha = k_{\nu_\alpha}$ which is invariant on $W$-orbits.  We also use the notation $X_b(q^a) :=  q^{(b,a)}$, and in particular, $X_b(q^{\rho_k}) = q^{(b,\rho_k)} = t^{(b,\rho_{\text{sht}})}_{\text{sht}}t^{(b,\rho_{\text{lng}})}_{\text{lng}}$.  Following \cite{C3}, we have the duality and evaluation formulas:
\begin{gather}
P_b(q^{c-\rho_k})P_c(q^{-\rho_k}) = P_c(q^{b-\rho_k})P_b(q^{-\rho_k})\hspace{5pt}\text{ for }\hspace{5pt}b,c\in P_-,\\
P_b(q^{-\rho_k}) = q^{-(\rho_k,b)}\displaystyle\prod_{\alpha\in R_+}\prod_{j=0}^{(\alpha^\vee,b)-1}\left(\frac{1-q_\alpha^jt_\alpha X_\alpha(q^{\rho_k})}{1-q_\alpha^j X_\alpha(q^{\rho_k})}\right).
\end{gather}
The corresponding \emph{spherical polynomial} is $P^\circ_b :=  P_b/P_b(q^{-\rho_k})$.

\subsection{DAHA-Jones polynomials}

Here we provide an efficient definition of the DAHA-Jones polynomials, which were originally defined in \cite{ChJ,CJJ} for torus knots and extended to iterated torus knots in \cite{CD}.  We also state their main (algebraic) properties, which were conjectured in \cite{ChJ} and mostly proved in \cite{CJJ,GN}.

\subsubsection{\sf $PSL^\wedge_2(\mathbb{Z})$-action}

Define a central idempotent:
\begin{equation}
e :=  \frac{1}{|W|}\displaystyle\sum_{w\in W}w, 
\end{equation}
in the group algebra of $W$.  Then the \emph{spherical DAHA} is $\mathcal{S\!H}:= e\mathcal{H\!\!\!H}e\subset\mathcal{H\!\!\!H}$.  In particular, $P^\circ_b\in\mathcal{S\!H}$.  Further, define the \emph{projective $PSL_2(\mathbb{Z})$} by
\begin{equation}
PSL^\wedge_2(\mathbb{Z}):= \langle\tau_\pm\text{ : } \tau_+\tau^{-1}_-\tau_+ = \tau^{-1}_-\tau_+\tau^{-1}_-\rangle,
\end{equation}
as a group whose action $\mathcal{H\!\!\!H}$ is represented by:
\begin{equation}\label{PSLaction}
\tau_+ = \left( \begin{array}{cc}
1 & 1\\
0 & 1\end{array} \right)\text{, }\tau_- = \left( \begin{array}{cc}
1 & 0\\
1 & 1\end{array} \right)\text{, where }\left( \begin{array}{cc}
a & b\\
c & d\end{array} \right): \begin{cases} 
      X_\lambda\mapsto X_\lambda^aY_\lambda^c\\
      T_i\mapsto T_i\\
      Y_\lambda\mapsto X_\lambda^bY_\lambda^d
   \end{cases},
\end{equation}
for $\lambda\in P$, $i>0$ and extends to an action on all of $\mathcal{H\!\!\!H}$, which restricts to an action on $\mathcal{S\!H}\subset\mathcal{H\!\!\!H}$.

\subsubsection{\sf Evaluation coinvariant}

We define a functional $\{\cdot\}_{ev}:\mathcal{H\!\!\!H}\rightarrow\mathbb{Z}_{q,t}$ called the \emph{evaluation coinvariant} which first writes $h\in\mathcal{H\!\!\!H}$,
\begin{equation}
h= \displaystyle\sum_{a,w,b}c_{a,w,b}X_a T_w Y_b,
\end{equation}
in the unique form guaranteed by the PBW Theorem \ref{PBW} and then substitutes
\begin{equation}\label{evalsub}
X_a\mapsto q^{-(\rho_k,a)}\text{, }\hspace{10pt}T_i\mapsto t_i^{\frac{1}{2}}\text{, }\hspace{10pt}Y_b\mapsto q^{(\rho_k,b)}.
\end{equation}
This process factors through the polynomial representation, which allows one to avoid making direct use of the PBW theorem (which can be rather complicated to implement).  In other words, $\{\cdot\}_{ev}$ is equivalent to projection onto the polynomial representation followed by the substitution \eqref{evalsub}.  See \cite{C102}.

\subsubsection{\sf Main definition}

Corresponding to the $(r,s)$-torus knot, choose an element $\hat{\gamma}_{r,s}\in PSL^\wedge_2(\mathbb{Z})$ which is \emph{any} word in $\tau_\pm$ that can be represented by
\begin{equation}
\gamma_{r,s} = \left( \begin{array}{cc}
r & \ast\\
s & \ast\end{array} \right),
\end{equation}
where the $\ast$ entries do not matter, since $\hat{\gamma}_{r,s}$ will act on a polynomial in $X_i$, see \eqref{PSLaction}.  For any root system $R$ and dominant weight $b\in P_+$, let
\begin{gather}\label{JD}
J\!D^R_{r,s}(b; q,t) := \{\hat{\gamma}_{r,s}(P_b)/P_b(q^{-\rho_k})\}_{ev},\\
\widetilde{J\!D}^R_{r,s}(b; q,t) :=  q^\bullet t^\bullet J\!D^R_{r,s}(b; q,t),
\end{gather}
where $q^\bullet t^\bullet$ is the lowest $q,t$-monomial in $J\!D^R_{r,s}(b; q,t)$, if it is well-defined.  Then $\widetilde{J\!D}^R_{r,s}(b; q,t)\in\mathbb{Z}[q,t]$ is the (reduced, tilde-normalized) \emph{DAHA-Jones polynomial}.

\subsubsection{\sf Properties of DAHA-Jones polynomials}

Here we recall some important properties of DAHA-Jones polynomials, which were conjectured in \cite{ChJ} and proved in Theorem 1.2 of \cite{CJJ}.  First, we remark that the tilde-normalized DAHA-Jones polynomials are, in fact, polynomials: 
\begin{equation}
\widetilde{J\!D}^R_{r,s}(b; q,t)\in\mathbb{Z}[q,t].
\end{equation}
Then, in anticipation of a connection to quantum knot invariants, we expect that DAHA-Jones polynomials should satisfy the usual topological properties with respect to the torus knot $T^{r,s}$:
\begin{enumerate}
\item \emph{(well-defined)} $\widetilde{J\!D}^R_{r,s}(b; q,t)$ does not depend on the choice of $\hat{\gamma}_{r,s}\in PSL^\wedge_2(\mathbb{Z})$,
\item \emph{(unknot)} $\widetilde{J\!D}^R_{r,1}(b; q,t) = 1$,
\item \emph{($r,s$-symmetry)} $\widetilde{J\!D}^R_{r,s}(b; q,t) = \widetilde{J\!D}^R_{s,r}(b; q,t)$,
\item \emph{(orientation)} $\widetilde{J\!D}^R_{r,s}(b; q,t) = \widetilde{J\!D}^R_{-r,-s}(b; q,t)$,
\item \emph{(mirror image)} $J\!D^R_{r,-s}(b; q,t) = J\!D^R_{r,s}(b; q^{-1},t^{-1})$.
\end{enumerate}
Finally, the following \emph{evaluation} is a property of the refinement which reflects ``exponential growth" in the number of terms in $\widetilde{J\!D}^R_{r,s}(b; q,t)$ with respect to $|b|$:
\begin{equation}\label{JDeval}
\widetilde{J\!D}^R_{r,s}(\displaystyle\sum_{i=1}^n b_i\omega_i; q=1,t) = \prod_{i=1}^n\widetilde{J\!D}^R_{r,s}(\omega_i; q=1,t)^{b_i}.
\end{equation}
It is related to the fact that $P_{b+c}=P_bP_c$ upon $q\rightarrow 1$.  We do not discuss the \emph{color exchange}, which is also part of Theorem 1.2 and corresponds to generalized level-rank duality.

\subsection{Relation to torus knot polynomials/homologies}
\subsubsection{\sf Quantum groups}
In \cite{ChJ} it was demonstrated for $A_1$, announced for $A_n$, and conjectured for general root systems that, upon $t\mapsto q$, the DAHA-jones polynomials coincide (up to an overall factor) with the corresponding (normalized/reduced) quantum invariants of torus knots:
\begin{equation}
q^{\bullet}\widetilde{J\!D}^{R}_{r,s}(b;q,q) \stackrel{\text{conj.}}{=\!=} P^{\mathfrak{g},V_b}(T^{r,s};q).
\end{equation}
Here $\mathfrak{g}$ is the Lie algebra corresponding to the root system $R$, and $V_b$ is the representation of $\mathfrak{g}$ with highest weight $b\in P_+(R)$.

In the author's (R.E.) Ph.D. thesis, this connection has been established for $R$ of types $A$ and $D$, as well as for the examples of $(E_6,\omega_1)$ used in this paper.

\subsubsection{\sf DAHA-superpolynomials}

Here we restrict to type-$A$ root systems and present the ``three super-conjectures" from Section 2.2 of \cite{ChJ}, which are now theorems due to \cite{CJJ,GN}.

\begin{thm}\label{DAHAsupthm} For any $n\geq m-1$, we may naturally interpret $\lambda\in P_+(A_m)$ as a weight for $A_n$.  

\begin{enumerate}
\item {\bf (Stabilization)} There exists a unique polynomial $H\!D_{r,s}(\lambda;q,t,a)\in\mathbb{Z}[q,t^{\pm 1},a]$, defined by the (infinitely many) specializations
\begin{equation}
H\!D_{r,s}(\lambda;q,t,a\mapsto -t^{n+1}) = \widetilde{J\!D}^{A_{n}}_{r,s}(\lambda;q,t)\text{, for }n\geq m-1.
\end{equation}
We will call $H\!D_{r,s}(\lambda;q,t,a)$ the \underline{DAHA-superpolynomial}.

\item {\bf (Duality)} Let $q^At^B$ be the greatest $q,t$-monomial in $H\!D_{r,s}(\lambda;q,t,a)$ whose $a$-degree is $0$.  Then
\begin{equation}
H\!D_{r,s}(\lambda^{tr};q,t,a) = t^Aq^BH\!D_{r,s}(\lambda;t^{-1},q^{-1},a),
\end{equation}
where $\lambda^{tr}$ indicates the transposed Young diagram for $\lambda$.

\item {\bf (Evaluation)} It immediately follows from \eqref{JDeval} that
\begin{equation}\label{qeval}
H\!D_{r,s}\Bigl(\displaystyle\sum_{i=1}^m \lambda_i\omega_i;1,t,a\Bigr) = \prod_{i=1}^m \Bigl(H\!D_{r,s}(\omega_i;1,t,a)\Bigr)^{\lambda_i}.
\end{equation}
When combined with the duality, this implies
\begin{equation}\label{teval}
H\!D_{r,s}\Bigl(\displaystyle\sum_{i=1}^m \lambda_i\omega_i;q,1,a\Bigr) = \prod_{i=1}^m \Bigl(H\!D_{r,s}(\omega_i;q,1,a)\Bigr)^{\lambda_i}.
\end{equation}
Currently, the latter has no direct interpretation in terms of Macdonald polynomials or the DAHA-Jones construction.
\end{enumerate}
\end{thm}

We can generally make contact with the conventions used in the literature on superpolynomials, e.g., \cite{DGR}, by a transformation $\text{DAHA}\mapsto\text{DGR}$:
\begin{equation}\label{DAHAtoDGR}
t \mapsto q^{2} \,, \qquad
q \mapsto q^2 t^2 \,, \qquad
a \mapsto a^2t.
\end{equation}
Then we have the following conjecture, which extends the one from \cite{AS}.

\begin{conj}\label{DAHADGR}
For a \emph{rectangular} Young diagram $i\times j$, i.e., a weight $j\omega_i\in P_+$, the coefficients of $H\!D_{r,s}(j\omega_i;q,t,a)$ are positive integers.  In this case, upon the transformation \eqref{DAHAtoDGR}, one recovers the superpolynomials from \cite{DGR,GS,GGS}.
\end{conj}

In light of Conjecture \ref{DAHADGR}, one can attribute the duality to the ``mirror symmetry" and the evaluation to the ``refined exponential growth" of \cite{GS,GGS}.  Furthermore, in Lemma 2.8 of \cite{GN} the authors demonstrate that the DAHA-Jones polynomials are proper (formal) generalizations--to any root system and weight--of the refined torus knot invariants of \cite{AS}.
\medskip

\setcounter{equation}{0}
\section{Exceptional knot homology}
\label{sec:hyper}

\subsection{Approach: DAHA + BPS}\label{approach}

In \cite{DGR} the authors introduce the superpolynomial for knot homologies,
as a generating function of the refined BPS invariants on the one hand and as
the Poincar\'e polynomial of the HOMFLY homology on the other.  Analogous
constructions for colored HOMFLY and Kauffman homologies were developed in
\cite{GS} and \cite{GW}, respectively.  Here, we incorporate the exceptional
Lie algebra $\mathfrak{e}_6$ and its $27$-dimensional representation with
(minuscule) highest weight $\omega_1$.

Exceptional Lie algebras pose a number of unique challenges.  For one, they are
singular in the sense that they do not belong to infinite families in any obvious
way.  Thus, we are missing a natural notion of ``stabilization," which helps the
identification of gradings/differentials in the classical cases.

In \cite{CE}, the authors consider stabilization for the
Deligne-Gross ``exceptional series."  However, this is a fundamentally different
phenomenon than considered here, as their examples contain negative coefficients.
It is an interesting question, relegated to future research, whether the approach
in \cite{CE} is compatible with the approach here.

We also face a more technical/computatational challenge.  Even the
ordinary (quantum group) knot invariants for $\mathfrak{e}_6$ have not been
explicitly computed in the literature.  The author R.E. has computed them for the
cases considered here (unpublished) and verified their coincidence with the DAHA-Jones
polynomials upon $t\mapsto q$. Furthermore, no corresponding homology theory has
been formally defined.

We manage to overcome these obstacles by applying the technique of differentials
from \cite{DGR, R} to the DAHA-Jones polynomials, $q,t$-counterparts of quantum knot
invariants defined in \cite{ChJ}.  This combination is sufficiently powerful to overcome
all obstacles.  Here, we propose new invariants, the \emph{hyperpolynomials}, for
$\mathfrak{e_6}, \mathbf{27}$ torus knot homologies, as well as produce some explicit examples.

\subsubsection{\sf Notation and conventions}
We will use two sets of conventions in this paper: the standard DAHA conventions and
conventions used in the literature on quantum group invariants (``QG").
While our calculations are performed in DAHA conventions ($q,t,a$), we are ultimately interested in QG conventions ($q,t,u$).  To change $\text{DAHA}\rightarrow\text{QG}$, we apply the ``grading change" isomorphism:
\begin{equation}\label{change}
a \mapsto u t^{-1} \,, \qquad
q \mapsto q t^2 \,, \qquad
t \mapsto q.
\end{equation}
Even though $q,t$ are used in both sets of conventions, whether we are referring to DAHA or QG will be contextually clear.

Furthermore, for a given knot, polynomials in QG conventions are usually associated to a Lie algebra $\mathfrak{g}$ and a representation ($\mathfrak{g}$-module) $V$.  Polynomials in DAHA conventions are (equivalently) associated to a root system $R$ and a (dominant) weight $b\in P_+$.  The correspondence between $\mathfrak{g}$ and $R$ is via the classification of complex, semisimple Lie algebras, and $b$ is the highest weight for $V$, as labeled in \cite{Bo}.
\smallskip

Now, in QG-conventions, our hyperpolynomials are Poincar\'{e} polynomials for a (hypothetical) triply-graded vector space:
\begin{equation}\label{HQG}
H^{\frak e_6, \mathbf{27}} (K;q,t,u) \; := \; \sum_{i,j,k} q^i t^j u^k \text{dim} \mathcal{H}^{\frak e_6, \mathbf{27}}_{i,j,k} (K).
\end{equation}
The usual two-variable Poincar\'e polynomials \eqref{PHgR} are returned upon setting $u=1$:
\begin{equation}\label{qtPon}
\mathcal{P}^{\frak e_6, \mathbf{27}} (K;q,t) \; := \; H^{\frak e_6, \mathbf{27}} (K;q,t,1),
\end{equation}
and we have, upon taking the graded Euler characteristic with respect to $t$,
\begin{equation}\label{QGinv}
P^{\frak e_6, \mathbf{27}} (K;q) \; = \; \mathcal{P}^{\frak e_6, \mathbf{27}} (K;q,-1),
\end{equation}
i.e. these ``categorify" the quantum knot invariants \eqref{PgR} for $\frak e_6, \mathbf{27}$.
\smallskip

This story may be translated into DAHA conventions.  In light of (\ref{change}), we may also write the hyperpolynomials in DAHA conventions:
\begin{equation}\label{HDAHA}
H\!D^{E_6}_{r,s} (\omega_1;q,t,a) \; := \; \sum_{i,j,k} q^{\frac{j+k}{2}} t^{\frac{2i-j+k}{2}} a^k \text{dim} \mathcal{H}^{\frak e_6, \mathbf{27}}_{i,j,k} (T^{r,s}).
\end{equation}
for the \emph{same} vector space as in (\ref{HQG}).  Though we do not consider a DAHA analogue of $\mathcal{P}^{\frak e_6, \mathbf{27}}$ here, we may obtain the DAHA-Jones polynomial by taking the graded Euler characteristic with respect to $a$:
\begin{equation}\label{qtreduce}
\widetilde{J\!D}^{E_6}_{r,s} (\omega_1; q,t) \; = \; H\!D^{E_6}_{r,s} (\omega_1;q,t,-1).
\end{equation}
Recall that the DAHA-Jones polynomials are $t$-refinements of the QG knot invariants.  They are (conjecturally) related by setting $t\mapsto q$:
\begin{equation}\label{qreduce}
P^{\frak e_6, \mathbf{27}} (T^{r,s};q) \; = \; \widetilde{J\!D}^{E_6}_{r,s} (\omega_1; q,q).
\end{equation}
Thus, we come full circle and make contact with the QG conventions at the level of polynomials.
\smallskip

For the convenience of the reader, our conventions and notations are summarized in the following commutative diagram.

\begin{equation}\label{CDconv}
\begin{aligned}
{\normalsize
\xymatrixcolsep{3pc}\xymatrixrowsep{2pc}\xymatrix{
\fbox{\text{DAHA}}&H\!D \ar[rr]^{a=-1}_{(\ref{qtreduce})} \ar@{{<}-{>}}[dd]_{(\ref{change})} &&\widetilde{J\!D} \ar[dd]^{t\mapsto q}_{(\ref{qreduce})}\\
\mathcal{H}_{i,j,k}\ar[ur]^{(\ref{HDAHA})}\ar[dr]_{(\ref{HQG})}&&&\\
\fbox{\text{QG}}&H \ar[r]^{u=1}_{(\ref{qtPon})}& \mathcal{P} \ar[r]^{t=-1}_{(\ref{QGinv})}& P}}
\end{aligned}
\end{equation}
\smallskip

\subsubsection{\sf Torus knots}

Presently, our approach is confined to the torus knots and links for which the DAHA-Jones polynomials are defined.  The reason for this limitation is algebraic from the DAHA point of view.  Here we will shed some light on it geometrically and physically.

In the BPS framework, something special happens when $K = T^{r,s}$ is a torus knot. Then, the five-brane theory in \eqref{M51} has an extra $R$-symmetry $U(1)_R$ that acts on ${\bf S}^3$
leaving the knot $K = T^{r,s}$ and, hence, the Lagrangian $L_K \subset T^* {\bf S}^3$ invariant.
Following \cite{AS}, we denote the quantum number corresponding to this symmetry by $S_R$,
and also introduce the generating function, {\it cf.} \eqref{HspacesPh}:
\begin{equation}
H^{\frak g, V} (K;q,t,u) \; := \; \text{Tr}_{\mathcal{H}_{\text{BPS}}} \; q^{P} t^{F} u^{S_R} \,.
\end{equation}
that ``counts'' refined BPS states in the setup \eqref{M51}.

{}From the perspective of \cite{ORS}, which is related to the DAHA approach, this extra variable / grading comes from the symmetry of the algebraic curve,
\begin{equation}
x^r \; = \; y^s,
\end{equation}
whose intersection with a unit sphere in $\mathbb{R}^4 \cong \mathbb{C}^2$ defines a $(r,s)$ torus knot $T^{r,s}$.

In either case, the origin of the extra grading (resp. variable $u$) has nothing to do with the choice of
homology (Khovanov, colored HOMFLY, or other); it simply comes from a very special choice of the knot (link)
and exists only for torus knots and links.

As a result, what for a generic knot $K$ might be a doubly-graded homology $\mathcal{H}^{\frak g, V}_{i,j} (K)$
for torus knot becomes a triply-graded homology $\mathcal{H}^{\frak g, V}_{i,j,k} (K)$, with an extra $u$-grading.
Likewise, what normally would be a triply-graded (say, HOMFLY or Kauffman) homology, for a torus knot $K = T^{r,s}$
becomes a quadruply-graded homology $\mathcal{H}^{\frak g, V}_{i,j,k,\ell} (T^{r,s})$, c.f. \cite{GGS}.

\subsubsection{\sf Hyper-lift}
We wish to elevate the two-variable DAHA-Jones polynomial $\widetilde{J\!D}^{E_6}_{r,s}(\omega_1; q,t)$, which in general has both positive and negative coefficients, to a three-variable hyperpolynomial $H\!D^{r,s}_{E_6}(\omega_1; q,t,a)$ with only positive coefficients.

As in (\ref{HDAHA}), this ``upgraded" polynomial will be the Poincar\'e polynomial of a triply-graded vector space $\mathcal{H}^{\mathfrak{e}_6,\mathbf{27}}_{i,j,k}(T^{r,s})$, accounting for its positive coefficients.  As in (\ref{qtreduce}), it is related to $\widetilde{J\!D}^{E_6}_{r,s}$ by taking the graded Euler characteristic with respect to the $k$-grading (resp. variable $a$):
\begin{equation}
\widetilde{J\!D}^{E_6}_{r,s}(\omega_1; q,t) = H\!D^{r,s}_{E_6}(\omega_1; q,t,-1).
\end{equation}
Note that we are here constructing the polynomial $H\!D^{r,s}_{E_6}(\omega_1)$ whose constituent monomials encode the graded dimensions of the irreducible components of the vector space $\mathcal{H}^{E_6,r,s}_{i,j,k}$.  We are not constructing this vector space itself.

Of course, there will be many polynomials $H\!D^{r,s}_{E_6}(\omega_1)$ that satisfy only the aforementioned properties.  We will define ours intelligently so that it is uniquely determined and so that like the HOMFLY-PT (``superpolynomial") and Kauffman homologies --- which, respectively, unify $\mathfrak{sl}_N$ and $\mathfrak{so}_N$ invariants --- our ``hyperpolynomial" will unify the $(\mathfrak{e}_6,\mathbf{27})$-invariant with invariants associated to ``smaller" algebras and representations $(\mathfrak{g},V)$.

\subsubsection{\sf Differentials and specializations}

This unification with other $(\mathfrak{g},V)$-colored invariants is realized via certain (conjectural) spectral sequences on $\mathcal{H}^{\mathfrak{e}_6,\mathbf{27}}_\ast$ induced by deformations of the potential $W_{E_6,27}\leadsto W_{\mathfrak{g},V}$, which are studied in section \ref{sec:Adj}.  With the additional assumption that these spectral sequences converge on its second page, such deformations gives rise to differentials $d_{\mathfrak{g},V}$ such that the homology:
\begin{equation}
H_\ast(\mathcal{H}^{\mathfrak{e}_6,\mathbf{27}}_\ast, d_{\mathfrak{g},V}) \cong \mathcal{H}^{\mathfrak{g},V}_\ast.
\end{equation}

Practically speaking, suppose that such a differential $d_{\mathfrak{g},V}$ exists ($=d_{R,b}$ in DAHA conventions),  and that its $(q,t,a)$-degree is $(\alpha,\beta, \gamma)$.  Then each monomial term in $H\!D^{r,s}_{E_6}(\omega_1)$ will participate in exactly one of two types of direct summands in the chain complex $(\mathcal{H}^{\mathfrak{e}_6,\mathbf{27}}_\ast, d_{R,b})$:
\begin{equation}
0\stackrel{d}{\longrightarrow} q^it^ja^k\stackrel{d}{\longrightarrow}0,
\end{equation}
\begin{equation}
0\stackrel{d}{\longrightarrow} q^it^ja^k\stackrel{\cong}{\longrightarrow} q^{i+\alpha}t^{j+\beta}a^{k+\gamma}\stackrel{d}{\longrightarrow}0.
\end{equation}
Observe that we can re-express this as a decomposition:
\begin{equation}\label{diffspec}
H\!D^{r,s}_{E_6}(\omega_1) = \widetilde{H\!D}_{R}(b) + (1+q^\alpha t^\beta a^\gamma)\mathcal{Q}(q,t,a),
\end{equation}
where $\widetilde{H\!D}_{R}(b)$ is related to $\widetilde{J\!D}^{R}(b)$ by the specialization:
\begin{equation}
H\!D^{r,s}_{E_6}(\omega_1; a = -q^{-\frac{\alpha}{\gamma}}t^{-\frac{\beta}{\gamma}}) = \widetilde{H\!D}_R(b;a = -q^{-\frac{\alpha}{\gamma}}t^{-\frac{\beta}{\gamma}}) = \widetilde{J\!D}^{R}(b)
\end{equation}
which subsumes the differential $d_{R,b}$, realized by setting $(1+q^\alpha t^\beta a^\gamma) = 0$.  Note that since these polynomials always have integer exponents (corresponding to integer gradings of a vector space), we will always be able to define the $a$-grading in such a way that $\gamma$ divides $\alpha$ and $\beta$.

To restore the $a$-grading to $\widetilde{J\!D}^{E_6}(\omega_1)$, we play this game in reverse.  On the $q,t$-level, we have a decomposition:
\begin{equation}
\widetilde{J\!D}^{E_6}(\omega_1) = \widetilde{J\!D}^{R}(b) + (1\pm q^\alpha t^\beta)Q(q,t).
\end{equation}
Since many of the polynomials $\widetilde{J\!D}^{R}(b)$ are known, we can hope to use this structure to recover the $a$-gradings of specific generators as well as the $a$-degrees of the $d_{R,b}$.  If we can do this for sufficiently many $(R,b)$, we will obtain enough constraints (specializations) to uniquely define the (relative) $a$-grading in $H\!D^{r,s}_{E_6}(\omega_1)$.

\subsubsection{\sf Uniqueness}
Suppose that we have defined $H\!D$ by some (possibly infinite) set of differentials/specializations $S:=\{(R,b, \alpha,\beta,\gamma)\}$, each of the form (\ref{diffspec}) with the \emph{same} $\widetilde{H\!D}_R(b)$.  If two hyperpolynomials $H\!D_1$, $H\!D_2$ each satisfy all of the specializations $S$, then evidently $H\!D_1-H\!D_2\in I_S$, where
\begin{equation}
I_S:=\prod_S \left(1+q^\alpha t^\beta a^\gamma\right)
\end{equation}
is an ideal in $R:=\mathbb{Z}[[q,t,a]]$.  Then $H\!D$ corresponds to a unique coset $[H\!D]\in R/I_S$.

If $S$ is infinite, then we may choose a distinguished representative of $[H\!D]$, i.e. the only one with finitely many terms.  This is precisely the situation when considering superpolynomials and hyperpolynomials for the classical series of Lie algebras.

When $S$ is finite, there is also a distinguished representative.  Since $H\!D$ is required to have \emph{positive} coefficients, we may simply require that it is minimal in $[H\!D]$ with respect to that property, i.e. it has the minimum number of terms.

Indeed, suppose $H\!D_1\neq H\!D_2$ are minimal, and write $H\!D_1-H\!D_2 = F\cdot\prod_S (1+q^\alpha t^\beta a^\gamma)\in I_S$ for some $F\in R$.  Since the $H\!D_i$ both have positive coefficients, we may write $F=F_1-F_2$, where each $F_i$ has only positive coefficients.  Then clearly the monomials in $F_i\cdot\prod_S (1+q^\alpha t^\beta a^\gamma)$ are all monomials in $H\!D_i$, and since these belong to $I_S$, they cancel in every specialization in $S$.  Then
\begin{equation}
H\!D_i' := H\!D_i - F_i\cdot\prod_S (1+q^\alpha t^\beta a^\gamma)
\end{equation}
is a new polynomial with positive coefficients, fewer terms, and which satisfies all of the specializations $S$.  This contradicts the assumed minimality of $H\!D_i$.

Restricting ourselves to these distinguished representatives, the uniqueness of our $H\!D$ depends on the uniqueness of the $\widetilde{H\!D}_R(b)$ chosen simultaneously for $(R,b)\in S$.  As we will see below, this is manifest in all cases considered.

\subsection{$E_6$-hyperpolynomials}

In the standard knot theory (QG) conventions, our main proposal for $H^{\frak e_6, \mathbf{27}}$ is based on the following (finite) set of differentials/specializations:

\begin{equation}
    \begin{tabular}{| l | l | l |}
    \hline
    $\mathfrak{g},V$ & $H^{\frak e_6, \mathbf{27}}(u=1, t=?) = P^{\mathfrak{g},V}$ & deg$(d_{\mathfrak{g},V})$ \\ \hline \hline
    $\mathfrak{e}_6,\mathbf{27}$ & $-1$ & $(0,-1,1)$  \\ \hline
    $\mathfrak{d}_5, \mathbf{10}$ & $-q^{4}$ & $(4,-1,1)$ \\ \hline
    $\mathfrak{a}_6, \mathbf{7}$ & $-q^{5}$ & $(5,-1,1)$  \\ \hline
    canceling & $-q^{8}$ & $(8,-1,1)$  \\ \hline
    canceling & $-q^{-13}$ & $(13,1,1)$  \\
    \hline
    \end{tabular}
\end{equation}
which we will take as a definition for our hyperpolynomial.  By a ``canceling" differential, we mean that the corresponding homology is one dimensional.  In other words,  $\widetilde{H\!D}_{R}(b)$ in (\ref{diffspec}) --- as well as its variant in the QG conventions --- is a single monomial.

We construct three explicit examples, for $T^{3,2}$, $T^{5,2}$, and $T^{4,3}$ torus knots, which are also known as the ${\bf 3}_1$, ${\bf 5}_1$, and ${\bf 8}_{19}$ knots, respectively. The result looks as follows:

\begin{equation}
H^{\frak e_6, \mathbf{27}}({\bf 3}_1) = 1 + q^2t^2 + q^5t^2 + q^{10}tu + q^{13}tu + q^{10}t^4 + q^{15}t^{3}u + q^{18}t^{3}u + q^{23}t^{2}u^2
\end{equation}

\begin{equation}
H^{\frak e_6, \mathbf{27}}({\bf 5}_1)\ =
\end{equation}
\renewcommand{\baselinestretch}{0.5}
{\small
\noindent\(
1 + q^2t^2+ q^5t^2 + q^{10}tu + q^{13}tu + q^4t^4 + q^7t^4 + q^{10}t^4
 + q^{12}t^{3}u + 2q^{15}t^{3}u + q^{18}t^{3}u + q^{23}t^{2}u^2 + q^{12}t^6 +q^{15}t^{6}
+ q^{17}t^{5}u + 2q^{20}t^{5}u +q^{23}t^{5}u + q^{25}t^{4}u^2 + q^{28}t^{4}u^2 + q^{20}t^{8}
 + q^{25}t^{7}u +q^{28}t^{7}u + q^{33}t^{6}u^2,
\)}
\renewcommand{\baselinestretch}{1.2}
\medskip

\begin{equation}
H^{\frak e_6, \mathbf{27}}({\bf 8}_{19})\ =
\end{equation}
\renewcommand{\baselinestretch}{0.5}
{\small
\noindent\(
1 + q^2t^2 + q^5t^2 + q^{10}tu + q^{13}tu + q^3t^4 + q^4t^4 + q^6t^4
 + q^7t^4 + q^{10}t^4 + q^{11}t^3u + q^{12}t^{3}u + q^{14}t^{3}u + 2q^{15}t^{3}u
 + q^{18}t^{3}u + q^{23}t^{2}u^2 + q^6t^6 + q^8t^6 + q^9t^6 + q^{11}t^6 + q^{12}t^6
 + q^{13}t^{5}u + q^{14}t^{5}u + q^{15}t^{6} + 3q^{16}t^{5}u + 2q^{17}t^{5}u + 2q^{19}t^{5}u
 + 2q^{20}t^{5}u + q^{21}t^{4}u^2 + q^{23}t^{5}u + 2q^{24}t^{4}u^2 + q^{25}t^{4}u^2
 + q^{27}t^{4}u^2 + q^{28}t^{4}u^2 + q^{12}t^8 + q^{13}t^8 + q^{14}t^{8} + q^{16}t^{8}
 + q^{17}t^{8} + q^{17}t^{7}u + q^{18}t^{7}u + q^{19}t^{7}u + q^{20}t^{8} + q^{20}t^{7}u
 + 3q^{21}t^{7}u + 2q^{22}t^{7}u + 2q^{24}t^{7}u + 2q^{25}t^{7}u + q^{25}t^{6}u^2
 + 2q^{26}t^{6}u^2 + q^{27}t^{6}u^2 + q^{28}t^{7}u + 3q^{29}t^{6}u^2 + q^{30}t^{6}u^2
 + q^{32}t^{6}u^2 + q^{33}t^{6}u^2 + q^{34}t^{5}u^3 + q^{37}t^{5}u^3 + q^{18}t^{10} + q^{21}t^{10}
 + q^{22}t^{10} + q^{22}t^{9}u + q^{23}t^{9}u + q^{25}t^{10} + q^{25}t^{9}u + 3q^{26}t^{9}u
 + q^{27}t^{9}u + q^{27}t^{8}u^2 + 2q^{29}t^{9}u + 2q^{30}t^{9}u + 2q^{30}t^{8}u^2
 + 2q^{31}t^{8}u^2 + q^{33}t^{9}u + q^{33}t^{8}u^2 + 3q^{34}t^{8}u^2 + q^{35}t^{8}u^2
 + q^{35}t^{7}u^3 + q^{37}t^{8}u^2 + q^{38}t^{8}u^2 + q^{38}t^{7}u^3 + q^{39}t^{7}u^3
 + q^{42}t^{7}u^3 + q^{30}t^{12} + q^{31}t^{11}u + q^{34}t^{11}u + q^{35}t^{11}u
 + q^{35}t^{10}u^2 + q^{36}t^{10}u^2 + q^{38}t^{11}u + 2q^{39}t^{10}u^2 + q^{40}t^{9}u^3
 + q^{42}t^{10}u^2 + q^{43}t^{10}u^2 + q^{43}t^{9}u^3 + q^{44}t^{9}u^3 + q^{47}t^{9}u^3
 + q^{48}t^{8}u^4.
\)}
\renewcommand{\baselinestretch}{1.2}
\medskip

Spectral sequence diagrams, which reveal the structure of the proposed differentials, are included for these examples in Appendix \ref{fig}.

\subsection{Computations with DAHA-Jones polynomials}

Here we demonstrate explicitly how the DAHA-Jones polynomials are combined with the theory of differentials to produce our examples.  First, we rewrite our proposal in DAHA conventions:
\begin{equation}
    \begin{tabular}{| l | l | l |}
    \hline
    $R,b$ & $H\!D^{E_6}_{r,s} (\omega_1;a=?) = \widetilde{J\!D}^{R}(b)$ & deg$(d_{R,b})$ \\ \hline \hline
    $E_6,\omega_1$ & $-1$ & $(0,0,1)$  \\ \hline
    $D_5, \omega_1$ & $-t^{-4}$ & $(0,4,1)$ \\ \hline
    $A_6, \omega_1$ & $-t^{-5}$ & $(0,5,1)$  \\ \hline
    canceling & $-t^{-8}$ & $(0,8,1)$  \\ \hline
    canceling & $-q^{-1}t^{-12}$ & $(1,12,1)$  \\
    \hline
    \end{tabular}
\end{equation}
This is identically our proposal for $H^{\frak e_6, \mathbf{27}}$ before the transformation (\ref{change}).  Now we consider each of our three examples individually.

\subsubsection{The Trefoil $T_{3,2}$}
The DAHA-Jones $(E_6,\omega_1)$ polynomial for the trefoil is
\begin{equation}
\widetilde{J\!D}^{E_6}_{3,2}(\omega_1; q,t) = 1 + qt + qt^4 - qt^9 - qt^{12} + q^2t^8 - q^2t^{13} - q^2t^{16} + q^2t^{21}.
\end{equation}
To elevate this to a Poincar\'e polynomial with positive coefficients, we introduce an extra $a$-grading.  For now this will only be a $\mathbb{Z}/2\mathbb{Z}$-grading ($a^{\underline{0}}$ or $a^{\underline{1}}$) compatible with the specialization $a=-1$:
\begin{equation}
\underline{H\!D}^{E_6}_{3,2} (\omega_1) = a^{\underline{0}} + qta^{\underline{0}} + qt^4a^{\underline{0}} + qt^9a^{\underline{1}} + qt^{12}a^{\underline{1}} + q^2t^8a^{\underline{0}} + q^2t^{13}a^{\underline{1}} + q^2t^{16}a^{\underline{1}} + q^2t^{21}a^{\underline{0}}.
\end{equation}

Now we would like to lift this $\mathbb{Z}/2\mathbb{Z}$-grading to a genuine $\mathbb{Z}$-grading, for which we use the differential structure outlined above.  Fortunately, this case is resolved rather easily by considering the $(D_n,\omega_1)$ DAHA-Jones polynomial:
\begin{equation}
\widetilde{J\!D}^{D_n}_{3,2}(\omega_1; q,t)  =  1 + qt + qt^{n-1} - qt^n - qt^{2n-2} + q^2t^{2n-2} - q^2t^{2n-1} - q^2t^{3n-3} + q^2t^{3n-2},
\end{equation}
which has the same dimension as $\widetilde{J\!D}^{E_6}_{3,2}$, so we can completely restore the $a$-grading by understanding just a single differential to some $(D_n, \omega_1)$, if one exists.

Indeed, such a differential to $(D_5, \omega_1)$ is indicated by the expression:
\begin{equation}\underline{H\!D}^{E_6}_{3,2} (\omega_1) \, =\label{T23D5} \end{equation} 
\noindent
\smallskip
\renewcommand{\baselinestretch}{0.5} 
{\normalsize
\(
   a^{\underline{0}} + qta^{\underline{0}} + qt^4a^{\underline{0}} + qt^5a^{\underline{1}} + qt^{8}a^{\underline{1}} + q^2t^8a^{\underline{0}} + q^2t^{9}a^{\underline{1}}
+ q^2t^{12}a^{\underline{1}} + q^2t^{13}a^{\underline{0}} + (1+t^4a^{\underline{1}})(qt^5a^{\underline{0}} + qt^{8}a^{\underline{0}} + q^2t^9a^{\underline{0}}
+ q^2t^{12}a^{\underline{0}} + q^2t^{13}a^{\underline{1}} + q^2t^{17}a^{\underline{1}}).
\)
}
\renewcommand{\baselinestretch}{1.2} 
\bigskip

\noindent Observe that the $a$-grading of this differential must be $1$ if the corresponding specialization is to contain only integer powers of $t$.  Thus, the $a$-grading of a generator corresponds to the number of canceling pairs of terms required to fit that generator into the expression above.  For example, the generator $qt^9a^{\underline{1}}$ is realized in \eqref{T23D5} as:
\begin{equation}
qt^9a^{\underline{1}} = qt^{5}a^{\underline{1}} + (1+t^4a^{\underline{1}})qt^5a^{\underline{0}},
\end{equation}
so its $a$-grading is $1$.  However, the generator $q^2t^{21}a^{\underline{0}}$ is realized in \eqref{T23D5} as:
\begin{equation}
q^2t^{21}a^{\underline{0}} = q^2t^{13}a^{\underline{0}} + (1+t^4a^{\underline{1}})(q^2t^{13}a^{\underline{1}} + q^2t^{17}a^{\underline{1}}),
\end{equation}
so its $a$-grading is $2$.  Overall, we restore the $a$-grading as a $\mathbb{Z}$-grading:
\begin{equation}
H\!D^{E_6}_{3,2} (\omega_1) = 1 + qt + qt^4 + qt^9a + qt^{12}a + q^2t^8 + q^2t^{13}a + q^2t^{16}a + q^2t^{21}a^2.
\end{equation}
Observe that, as desired, we so far have the following specializations which determine the $a$-grading:
\begin{gather}
H\!D^{E_6}_{3,2} (\omega_1;a=-1) = \widetilde{J\!D}^{E_6}_{3,2}(\omega_1),\\
H\!D^{E_6}_{3,2} (\omega_1;a=-t^{-4}) = \widetilde{J\!D}^{D_5}_{3,2}(\omega_1).
\end{gather}
We also find the two canceling differentials:
\begin{gather}
H\!D^{E_6}_{3,2} (\omega_1;a=-t^{-8}) = 1,\\
H\!D^{E_6}_{3,2} (\omega_1;a=-q^{-1}t^{-12}) = q^2t^8,
\end{gather}
as well as the differential to $(A_6,\omega_1)$:
\begin{equation}
H\!D^{E_6}_{3,2} (\omega_1;a=-t^{-5}) = \widetilde{J\!D}^{A_6}_{3,2}(\omega_1).
\end{equation}

\subsubsection{The Torus Knot $T_{5,2}$}
We repeat the above construction for $T^{5,2}$ and restore the $a$-grading to $\widetilde{J\!D}^{E_6}_{5,2}(\omega_1)$ in a way that includes all of the same structure.  We have the DAHA-Jones $(E_6,\omega_1)$ polynomial for $T^{5,2}$:
\begin{equation}\widetilde{J\!D}^{E_6}_{5,2}(\omega_1; q,t)  \; =\end{equation} 
\noindent
\smallskip
\renewcommand{\baselinestretch}{0.5} 
{\small
\(
   1 + qt + qt^4 - qt^9 - qt^{12} + q^2t^2 + q^2t^{5} + q^2t^{8} - q^2t^{10}
- 2q^2t^{13} - q^2t^{16} + q^2t^{21} + q^3t^9 + q^3t^{12} - q^3t^{14} - 2q^3t^{17}
- q^3t^{20} + q^3t^{22} + q^3t^{25} + q^4t^{16} - q^4t^{21} - q^4t^{24} + q^4t^{29}.
\)
}
\renewcommand{\baselinestretch}{1.2} 
\bigskip

\noindent As above, we introduce a mod-2 grading compatible with the specialization $a=-1$:
\begin{equation}\underline{H\!D}^{E_6}_{5,2} (\omega_1)  \;=\end{equation} 
\noindent
\smallskip
\renewcommand{\baselinestretch}{0.5} 
{\small
\(
   a^{\underline{0}} + qta^{\underline{0}} + qt^4a^{\underline{0}} + qt^9a^{\underline{1}} + qt^{12}a^{\underline{1}} + q^2t^2a^{\underline{0}} + q^2t^{5}a^{\underline{0}} + q^2t^{8}a^{\underline{0}} 
  + q^2t^{10}a^{\underline{1}} + 2q^2t^{13}a^{\underline{1}} + q^2t^{16}a^{\underline{1}} + q^2t^{21}a^{\underline{0}} + q^3t^9a^{\underline{0}} + q^3t^{12}a^{\underline{0}}
  + q^3t^{14}a^{\underline{1}} + 2q^3t^{17}a^{\underline{1}} + q^3t^{20}a^{\underline{1}} + q^3t^{22}a^{\underline{0}} + q^3t^{25}a^{\underline{0}} + q^4t^{16}a^{\underline{0}}
  + q^4t^{21}a^{\underline{1}} + q^4t^{24}a^{\underline{1}} + q^4t^{29}a^{\underline{0}}.
\)
}
\renewcommand{\baselinestretch}{1.2} 
\bigskip

\noindent The $D_5$ DAHA-Jones is:
\begin{equation}\widetilde{J\!D}^{D_5}_{5,2}(\omega_1; q,t)  \;=\end{equation} 
\noindent
\smallskip
\renewcommand{\baselinestretch}{0.5} 
{\normalsize
\(
   1 + qt + qt^4 - qt^5 - qt^8 + q^2t^2 + q^2t^5 - q^2t^6 + q^2t^8
  - 2q^2t^9 - q^2t^{12} + q^2t^{13} + q^3t^9 - q^3t^{10} + q^3t^{12} - 2q^3t^{13}
  + q^3t^{14} - q^3t^{16} + q^3t^{17} + q^4t^{16} - q^4t^{17} - q^4t^{20} + q^4t^{21},
\)
}
\renewcommand{\baselinestretch}{1.2} 
\bigskip

\noindent which again has the same dimension as $\widetilde{J\!D}^{E_6}_{5,2}$, so we can restore the $a$-grading in the same manner:
\begin{equation}H\!D^{E_6}_{5,2} (\omega_1) \; =\end{equation} 
\noindent
\smallskip
\renewcommand{\baselinestretch}{0.5} 
{\normalsize
\(
   1 + qt + qt^4 + qt^9a + qt^{12}a + q^2t^2 + q^2t^5 + q^2t^8 + q^2t^{10}a + 2q^2t^{13}a  + q^2t^{16}a + q^2t^{21}a^2 + q^3t^9 +q^3t^{12} + q^3t^{14}a + 2q^3t^{17}a +q^3t^{20}a  + q^3t^{22}a^2 + q^3t^{25}a^2 + q^4t^{16} + q^4t^{21}a +q^4t^{24}a + q^4t^{29}a^2.
\)
}
\renewcommand{\baselinestretch}{1.2} 
\bigskip

\noindent Observe that, as with the trefoil, we have specializations:
\begin{gather}
H\!D^{E_6}_{5,2} (\omega_1;a=-1) = \widetilde{J\!D}^{E_6}_{5,2}(\omega_1),\\
H\!D^{E_6}_{5,2} (\omega_1;a=-t^{-4}) = \widetilde{J\!D}^{D_5}_{5,2}(\omega_1),\\
H\!D^{E_6}_{5,2} (\omega_1;a=-t^{-8}) = 1,\\
H\!D^{E_6}_{5,2} (\omega_1;a=-q^{-1}t^{-12}) = q^4t^{16},\\
H\!D^{E_6}_{5,2} (\omega_1;a=-t^{-5}) = \widetilde{J\!D}^{A_6}_{5,2}(\omega_1).
\end{gather}

\subsubsection{The Torus Knot $T_{4,3}$}

We have the DAHA-Jones $(E_6,\omega_1)$ polynomial for $T^{4,3}$:
\begin{equation}\widetilde{J\!D}^{E_6}_{4,3}(\omega_1; q,t)  =\end{equation} 
\noindent
\smallskip
\renewcommand{\baselinestretch}{0.5} 
{\small
\(
   1 + qt + qt^4 - qt^9 - qt^{12} + q^2t + q^2t^2 + q^2t^4 + q^2t^5 + q^2t^8 - q^2t^9 - q^2t^{10}  - q^2t^{12} - 2q^2t^{13} - q^2t^{16} + q^2t^{21} + q^3t^3 + q^3t^5 + q^3t^6 + q^3t^8 + q^3t^9 - q^3t^{10} - q^3t^{11} + q^3t^{12}  - 3q^3t^{13} - 2q^3t^{14} - 2q^3t^{16} - 2q^3t^{17} + q^3t^{18} - q^3t^{20} + 2q^3t^{21} + q^3t^{22} + q^3t^{24} + q^3t^{25} + q^4t^8  + q^4t^9 + q^4t^{10} + q^4t^{12} - q^4t^{14} - q^4t^{15} - 3q^4t^{17} - 2q^4t^{18} - 2q^4t^{20} - q^4t^{21} + 2q^4t^{22} + q^4t^{23}  - q^4t^{24} + 3q^4t^{25} + q^4t^{26} + q^4t^{28} +q^4t^{29} - q^4t^{30} - q^4t^{33} + q^5t^{13} + q^5t^{16} - q^5t^{18} - 3q^5t^{21}  - 2q^5t^{24} + 2q^5t^{26} + 3q^5t^{29} + q^5t^{32} - q^5t^{34} - q^5t^{37} + q^6t^{24} - q^6t^{25} - q^6t^{28} + q^6t^{30} - q^6t^{32}  + 2q^6t^{33} - q^6t^{34} + q^6t^{36} - q^6t^{38} - q^6t^{41} + q^6t^{42},
\)
}
\renewcommand{\baselinestretch}{1.2} 
\smallskip

\noindent and the $D_5$ DAHA-Jones is:
\begin{equation}\widetilde{J\!D}^{D_5}_{4,3}(\omega_1; q,t)  \;=\end{equation} 
\noindent
\smallskip
\renewcommand{\baselinestretch}{0.5} 
{\small
\(
   1 + qt + qt^4 - qt^5 - qt^8 + q^2t + q^2t^2 + q^2t^4 - q^2t^6
  - 2q^2t^9 - q^2t^{12} + q^2t^{13} + q^3t^3 + q^3t^5 - q^3t^7 + q^3t^8
  - 2q^3t^9 - q^3t^{10} - q^3t^{12} + q^3t^{14} + q^3t^{17} + q^4t^8 - q^4t^{11}
  - q^4t^{13} + q^4t^{15} - q^4t^{16} + q^4t^{17}.
\)
}
\renewcommand{\baselinestretch}{1.2} 
\bigskip

\noindent From the outset it is apparent that these do not have the same dimension, so the same approach will be less effective.  However, we can try to assign a monomial in $\widetilde{J\!D}^{E_6}_{4,3}(\omega_1)$ to each monomial in $\widetilde{J\!D}^{D_5}_{4,3}(\omega_1)$ so that they coincide in the specialization $a=-t^{-4}$.  That is, we consider the following subset of $\underline{H\!D}_{E_6}(\omega_1)$:
\begin{equation}\underline{H\!D}_{D_5/E_6}  \;=\end{equation} 
\noindent
\smallskip
\renewcommand{\baselinestretch}{0.5} 
{\normalsize
\(
a^{\underline{0}} + qta^{\underline{0}} + qt^4a^{\underline{0}} + qt^9a^{\underline{1}} + qt^{12}a^{\underline{1}} + q^2ta^{\underline{0}} + q^2t^2a^{\underline{0}} + q^2t^4a^{\underline{0}}
 + q^2t^{10}a^{\underline{1}}  + 2q^2t^{13}a^{\underline{1}} + q^2t^{16}a^{\underline{1}} + q^2t^{21}a^{\underline{0}} + q^3t^3a^{\underline{0}} + q^3t^5a^{\underline{0}}
 + q^3t^{11}a^{\underline{1}} + q^3t^8a^{\underline{0}}  + 2q^3t^{13}a^{\underline{1}} + q^3t^{14}a^{\underline{1}} + q^3t^{16}a^{\underline{1}} + q^3t^{22}a^{\underline{0}}
 + q^3t^{25}a^{\underline{0}} + q^4t^8a^{\underline{0}} + q^4t^{15}a^{\underline{1}}  + q^4t^{17}a^{\underline{1}} + q^4t^{23}a^{\underline{0}} + q^4t^{20}a^{\underline{1}}
 + q^4t^{25}a^{\underline{0}},
\)
}
\renewcommand{\baselinestretch}{1.2} 
\bigskip

\noindent which should specialize to $\widetilde{J\!D}^{D_5}_{4,3}$, and thus lifts to:
\begin{equation}H\!D_{D_5/E_6} \;  =\end{equation} 
\noindent
\smallskip
\renewcommand{\baselinestretch}{0.5} 
{\small
\(
   1 + qt + qt^4 + qt^9a + qt^{12}a + q^2t + q^2t^2 + q^2t^4 + q^2t^{10}a  + 2q^2t^{13}a + q^2t^{16}a + q^2t^{21}a^2 + q^3t^3 + q^3t^5 + q^3t^{11}a + q^3t^8  + 2q^3t^{13}a + q^3t^{14}a + q^3t^{16}a + q^3t^{22}a^2 + q^3t^{25}a^2 + q^4t^8 + q^4t^{15}a  + q^4t^{17}a + q^4t^{23}a^2 + q^4t^{20}a + q^4t^{25}a^2.
\)
}
\renewcommand{\baselinestretch}{1.2} 
\bigskip

\noindent Now we turn our eye to the complementary subset:
\begin{equation}\underline{H\!D}_{E_6\backslash D_5}  \; =  \end{equation} 
\noindent
\smallskip
\renewcommand{\baselinestretch}{0.5} 
{\small
\(
  q^2t^5a^{\underline{0}} + q^2t^8a^{\underline{0}} + q^2t^9a^{\underline{1}} + q^2t^{12}a^{\underline{1}} + q^3t^6a^{\underline{0}} + q^3t^9a^{\underline{0}} + q^3t^{10}a^{\underline{1}} + q^3t^{12}a^{\underline{0}} + q^3t^{13}a^{\underline{1}} + q^3t^{14}a^{\underline{1}} + q^3t^{16}a^{\underline{1}} + 2q^3t^{17}a^{\underline{1}} + q^3t^{18}a^{\underline{0}} + q^3t^{20}a^{\underline{1}} + 2q^3t^{21}a^{\underline{0}} + q^3t^{24}a^{\underline{0}} + q^4t^9a^{\underline{0}} + q^4t^{10}a^{\underline{0}} + q^4t^{12}a^{\underline{0}} + q^4t^{14}a^{\underline{1}} + 2q^4t^{17}a^{\underline{1}} + 2q^4t^{18}a^{\underline{1}} + q^4t^{20}a^{\underline{1}} + q^4t^{21}a^{\underline{1}} + 2q^4t^{22}a^{\underline{0}} + q^4t^{24}a^{\underline{1}} + 2q^4t^{25}a^{\underline{0}} + q^4t^{26}a^{\underline{0}} + q^4t^{28}a^{\underline{0}} + q^4t^{29}a^{\underline{0}} + q^4t^{30}a^{\underline{1}} + q^4t^{33}a^{\underline{1}} + q^5t^{13}a^{\underline{0}} + q^5t^{16}a^{\underline{0}} + q^5t^{18}a^{\underline{1}} + 3q^5t^{21}a^{\underline{1}} + 2q^5t^{24}a^{\underline{1}} + 2q^5t^{26}a^{\underline{0}} + 3q^5t^{29}a^{\underline{0}} + q^5t^{32}a^{\underline{0}} + q^5t^{34}a^{\underline{1}} + q^5t^{37}a^{\underline{1}} + q^6t^{24}a^{\underline{0}} + q^6t^{25}a^{\underline{1}} + q^6t^{28}a^{\underline{1}} + q^6t^{30}a^{\underline{0}} + q^6t^{32}a^{\underline{1}} + 2q^6t^{33}a^{\underline{0}} + q^6t^{34}a^{\underline{1}} + q^6t^{36}a^{\underline{0}} + q^6t^{38}a^{\underline{1}} + q^6t^{41}a^{\underline{1}} + q^6t^{42}a^{\underline{0}}. 
\)
}
\renewcommand{\baselinestretch}{1.2} 
\bigskip

\noindent We can use the degrees of the differentials (now known) to restore the $a$-grading on these generators.  For example, $q^2t\in H\!D_{D_5/E_6}$ and $q^2t^9a^{\underline{1}}\in\underline{H\!D}_{E_6\backslash D_5}$ should cancel in the differential of degree $(0,8,1)$, so we restore the $a$-degree $q^2t^9a$ on that generator.  Carrying this out fully, we obtain:
\begin{equation}H\!D_{E_6\backslash D_5} \; = \end{equation} 
\noindent
\smallskip
\renewcommand{\baselinestretch}{0.5} 
{\small
\(
   q^2t^5 + q^2t^8 + q^2t^9a + q^2t^{12}a + q^3t^6 + q^3t^9 + q^3t^{10}a + q^3t^{12} + q^3t^{13}a + q^3t^{14}a  + q^3t^{16}a + 2q^3t^{17}a + q^3t^{18}a^2 + q^3t^{20}a + 2q^3t^{21}a^2 + q^3t^{24}a^2 + q^4t^9 + q^4t^{10} + q^4t^{12} + q^4t^{14}a  + 2q^4t^{17}a + 2q^4t^{18}a + q^4t^{20}a + q^4t^{21}a + 2q^4t^{22}a^2 + q^4t^{24}a + 2q^4t^{25}a^2 + q^4t^{26}a^2 + q^4t^{28}a^2  + q^4t^{29}a^2 + q^4t^{30}a^3 + q^4t^{33}a^3 + q^5t^{13} + q^5t^{16} + q^5t^{18}a + 3q^5t^{21}a + 2q^5t^{24}a + 2q^5t^{26}a^2  + 3q^5t^{29}a^2 + q^5t^{32}a^2 + q^5t^{34}a^3 + q^5t^{37}a^3 + q^6t^{24} + q^6t^{25}a + q^6t^{28}a + q^6t^{30}a^2 + q^6t^{32}a  + 2q^6t^{33}a^2 + q^6t^{34}a^3 + q^6t^{36}a^2 + q^6t^{38}a^3 + q^6t^{41}a^3 + q^6t^{42}a^4.
\)
}
\renewcommand{\baselinestretch}{1.2} 
\bigskip

\noindent Finally, observe that some generators that should cancel in certain specializations do not.  For example, $q^4t^9$ should cancel in the differential of degree $(0,4,1)$, but there is no $q^4t^{13}a$.  Taking all differentials into account, we add the generators: \begin{gather}
\{q^4t^{13}, q^4t^{13}a, q^4t^{16}a, q^4t^{21}a^2, q^5t^{17}a, q^4t^{16}, q^5t^{17}, 2q^5t^{25}a^2, q^5t^{20},\notag\\q^5t^{28}a^2, q^5t^{22}a^2, q^5t^{20}a, q^4t^{21}a, q^5t^{22}a, q^5t^{30}a^3, 2q^5t^{25}a, q^5t^{33}a^3,\label{extragen}\\q^6t^{29}a^2, q^5t^{28}aq^6t^{29}a, q^6t^{37}a^3, q^5t^{30}a^2, q^5t^{33}a^2, q^6t^{37}a^2\},\notag
\end{gather}
and take the sum $H\!D_{D_5/E_6} + H\!D_{E_6\backslash D_5} +$ \eqref{extragen} to obtain:
\begin{equation}H\!D^{E_6}_{4,3} (\omega_1) \; = \end{equation} 
\noindent
\smallskip
\renewcommand{\baselinestretch}{0.5} 
{\small
\(
   1 + qt + qt^4 + qt^9a + qt^{12}a + q^2t + q^2t^2 + q^2t^4 + q^2t^5 + q^2t^8 + q^2t^9a  + q^2t^{10}a + q^2t^{12}a + 2q^2t^{13}a + q^2t^{16}a + q^2t^{21}a^2 + q^3t^3 + q^3t^5 + q^3t^6 + q^3t^8 + q^3t^9  + q^3t^{10}a + q^3t^{11}a + q^3t^{12} + 3q^3t^{13}a + 2q^3t^{14}a + 2q^3t^{16}a + 2q^3t^{17}a + q^3t^{18}a^2  + q^3t^{20}a + 2q^3t^{21}a^2 + q^3t^{22}a^2 + q^3t^{24}a^2 + q^3t^{25}a^2 + q^4t^8 + q^4t^9 + q^4t^{10} + q^4t^{12}  + q^4t^{13} + q^4t^{13}a + q^4t^{14}a + q^4t^{15}a + q^4t^{16} + q^4t^{16}a + 3q^4t^{17}a + 2q^4t^{18}a + 2q^4t^{20}a  + 2q^4t^{21}a + q^4t^{21}a^2 + 2q^4t^{22}a^2 + q^4t^{23}a^2 + q^4t^{24}a + 3q^4t^{25}a^2 + q^4t^{26}a^2 + q^4t^{28}a^2  + q^4t^{29}a^2 + q^4t^{30}a^3 + q^4t^{33}a^3 + q^5t^{13} + q^5t^{16} + q^5t^{17} + q^5t^{17}a + q^5t^{18}a + q^5t^{20}  + q^5t^{20}a + 3q^5t^{21}a + q^5t^{22}a + q^5t^{22}a^2 + 2q^5t^{24}a + 2q^5t^{25}a + 2q^5t^{25}a^2 + 2q^5t^{26}a^2  + q^5t^{28}a + q^5t^{28}a^2 + 3q^5t^{29}a^2 + q^5t^{30}a^2 + q^5t^{30}a^3 + q^5t^{32}a^2 + q^5t^{33}a^2 + q^5t^{33}a^3  + q^5t^{34}a^3 + q^5t^{37}a^3 + q^6t^{24} + q^6t^{25}a + q^6t^{28}a + q^6t^{29}a + q^6t^{29}a^2 + q^6t^{30}a^2 + q^6t^{32}a  + 2q^6t^{33}a^2 + q^6t^{34}a^3 + q^6t^{36}a^2 + q^6t^{37}a^2 + q^6t^{37}a^3 + q^6t^{38}a^3 + q^6t^{41}a^3 + q^6t^{42}a^4,
\)
}
\renewcommand{\baselinestretch}{1.2} 
\bigskip

\noindent and verify that it satisfies:
\begin{gather}
H\!D^{E_6}_{4,3} (\omega_1;a=-1) = \widetilde{J\!D}^{E_6}_{4,3}(\omega_1),\\
H\!D^{E_6}_{4,3} (\omega_1;a=-t^{-4}) = \widetilde{J\!D}^{D_5}_{4,3}(\omega_1),\\
H\!D^{E_6}_{4,3} (\omega_1;a=-t^{-8}) = 1,\\
H\!D^{E_6}_{4,3} (\omega_1;a=-q^{-1}t^{-12}) = q^6t^{24},\\
H\!D^{E_6}_{4,3} (\omega_1;a=-t^{-5}) = \widetilde{J\!D}^{A_6}_{4,3}(\omega_1).
\end{gather}
\smallskip

\subsection{Further properties}
We observe that our hyperpolynomials exhibit a number of potentially meaningful structures beyond their defining specializations/differentials.  

\subsubsection{\sf Dimensions} First, observe that
\begin{equation}\label{dimensions}
H\!D^{E_6}_{r,s} (\omega_1;q,\pm 1,a) = H\!D^{A}_{r,s} (\omega_2;q,\pm 1,a)
\end{equation}
in all examples considered, in spite of the fact that the weight $\omega_2$ for $A_n$ is non-minuscule.  These relations generalize the special evaluations at $t=1$ of DAHA-Jones polynomials and DAHA-superpolynomials.  In particular, using the evaluation and super-duality theorems from \cite{ChJ}, equation (\ref{dimensions}) implies that
\begin{equation}
H\!D^{E_6}_{r,s} (\omega_1;q,1,a) = \left(H\!D^{A}_{r,s} (\omega_1;q,1,a)\right)^2.
\end{equation}
In turn, we see that the \emph{dimensions}
\begin{equation}
\dim H\!D^{E_6}_{r,s} := H\!D^{E_6}_{r,s}(\omega_1;1,1,1)
\end{equation}
are perfect squares.  The dimensions for our examples $T^{3,2},T^{5,2},T^{4,3}$ are $9, 25,121$, respectively.  These properties are analogues of the \emph{refined exponential growth} \cite{GS,GGS} for the exceptional groups.

\subsubsection{\sf Hat symmetry}
We also have a ``hat symmetry" corresponding to the involution of the Dynkin diagram for $E_6$ which sends $\omega_1\mapsto\omega_6$.  We define
\begin{equation}
\widehat{H\!D}_{r,s}^{E_6}(\omega_1; q,t,a) := H\!D_{r,s}^{E_6}(\omega_1; q\mapsto q t^{4}, t, a\mapsto at^{-4}),
\end{equation}
which satisfies the specializations
\begin{gather}
\widehat{H\!D}_{r,s}^{E_6}(\omega_1; q,t,-1) = \widetilde{J\!D}^{E_6}_{r,s}(\omega_6; q,t),\\
\widehat{H\!D}_{r,s}^{E_6}(\omega_1; q,t,-t^{-1}) =  \widetilde{J\!D}^{A_6}_{r,s}(\omega_1; qt^4,t),\\
\widehat{H\!D}_{r,s}^{E_6}(\omega_1; q,t,-t^{-4}) = 1,\\
\widehat{H\!D}_{r,s}^{E_6}(\omega_1; q,t,-q^{-1}t^{-12}) = q^\alpha t^\beta.
\end{gather}

\subsubsection{\sf Other evaluations}

We also have another potentially meaningful specialization of our hyperpolynomials at $a=q^{-1}t^{-9}$:
\begin{gather}
H\!D^{E_6}_{3,2} (\omega_1;a=-q^{-1}t^{-9}) = qt - q^2t^7 + q^2t^8,\\
H\!D^{E_6}_{5,2} (\omega_1;a=-q^{-1}t^{-9}) = q^2t^2 - q^2t^8 + q^3t^9 - q^3t^{15} + q^4t^{16},\\
H\!D^{E_6}_{4,3} (\omega_1;a=-q^{-1}t^{-9})  =  q^3t^3 - q^3t^9 + q^4t^{10} - q^4t^{16} + q^5t^{13} - q^5t^{19}  + q^5t^{17}\\ - q^5t^{23} - q^4t^{12} + q^6t^{24}.\notag
\end{gather}
We do not recognize the resulting polynomials.  However, observe the significant reduction in the number of terms, as well as their regularity.
\medskip

\setcounter{equation}{0}
\section{Adjacency tree of the corank-2 singularity $Z_{3,0}$}
\label{sec:Adj}

In the previous section, we encountered several ``exceptional'' differentials that relate homological invariants
of knots colored by representations of exceptional groups to knot homologies associated with classical groups.
In this section we explain the origin of such differentials.

There are two general ways to predict {\it a priori} the structure of the differentials, both of which are rooted in physics.
One approach \cite{GS} involves analysis of the spectrum \eqref{HspacesPh} of BPS states (a.k.a. $Q$-cohomology)
and how it changes when one varies stability parameters, such as the K\"ahler modulus \eqref{tVolN}.
The second approach \cite{Go} is based on deformations of the Landau-Ginzburg potential,
which for the $27$-dimensional representation of $\mathfrak{g} = \mathfrak{e}_6$ has the form \cite{GW}
\begin{equation}
W_{E_6, 27} \; = \; z_1^{13} - \frac{25}{169} z_1 z_4^3 + z_4 z_1^9.
\label{WE627}
\end{equation}
In general (and in every physics-based approach to knot homology), homology of the unknot can be represented
as a $Q$-cohomology, i.e., the space of $Q$-closed but not $Q$-exact states (called BPS states) in a two-dimensional
theory on a cylinder, $\mathbb{R} \times (\text{unknot}) = \mathbb{R} \times \mathbf{S}^1$.
In some cases, this two-dimensional theory admits a Landau-Ginzburg description, which for certain
Lie algebras $\mathfrak{g}$ and representations $V$ has been identified in \cite{GW}.
In this approach, spectral sequences and differentials correspond to relevant deformations and RG flows of
the two-dimensional ``unknot theory'' which, in the Landau-Ginzburg description, simply manifest
as deformations of the potential.

Therefore, in our present problem we need to explore deformations of the potential \eqref{WE627}
which correspond to the adjacencies of the singularity $Z_{3,0}$.  Additionally, we perform a nontrivial verification of our calculations using the adjacency
of the spectra of singularities.  A good general reference for material in this section is \cite{AGV}.

\subsection{Singularities and Adjacency}

A singularity is an analytic apparatus that captures the local geometry of a holomorphic (smooth) function at a critical point.  For our purposes, we will consider functions $f:\mathbb{C}^n\rightarrow \mathbb{C}$ and without loss of generality, critical points at $0\in\mathbb{C}^n$.
\smallskip

Let $\mathcal{O}_n$ be the space of all germs at $0\in\mathbb{C}^n$ of holomorphic functions $f: \mathbb{C}^n\rightarrow \mathbb{C}$.  Then the group of germs of diffeomorphisms (biholomorphic maps) $g:(\mathbb{C}^n,0)\rightarrow (\mathbb{C}^n,0)$ acts on $\mathcal{O}_n$ by $g\cdot f = f\circ g^{-1}$.  The orbits of this action define equivalence classes in $\mathcal{O}_n$, and those classes for which $0$ is a critical point are called \emph{singularities}.  Consider a class $L$ as a subspace of $\mathcal{O}_n$.  An $l$\emph{-parameter deformation} of $f \in L\subset \mathcal{O}_n$ with \emph{base} $\Lambda = \mathbb{C}^l$ is the germ of a smooth map $F: \Lambda \rightarrow \mathcal{O}_n$ such that $F(0)=f$.  
\smallskip

If $L$ is contained in the closure of some other subspace, $L\subset \bar{K}\subset \mathcal{O}_n$, then an infinitesimal neighborhood of every $f \in L\subset \mathcal{O}_n$ intersects $K$ nontrivially.  This geometric notion can be reformulated equivalently in terms of deformations and gives rise to the concept of adjacency.  That is, suppose that every function $f\in L$ can be transformed to a function in the class $K$ by an arbitrarily small deformation.  Here the ``size" of a deformation is a restriction on $\lambda\in\Lambda$, induced by the standard metric on $\mathbb{C}^l$.  In this case, we say that the singularity classes $L$, $K$ are \emph{adjacent}, written $L\rightarrow K$.

\subsubsection{\sf Versal deformations}
Here we aim to find the adjacencies to the specific class $Z_{3,0}$, that is the classes $K$ such that $Z_{3,0}\rightarrow K$.  We go about this by considering a specific type of deformation.

A deformation $F:\Lambda\rightarrow\mathcal{O}_n$ of $f$ is \emph{versal} if every deformation of $f$ is equivalent to one induced (by change of base $\Lambda$) from $F$.  If, in addition, $\Lambda$ has the smallest possible dimension, $F$ is said to be \emph{miniversal}, i.e., ``minimal and universal."

We can construct an explicit miniversal deformation of $f\in L$ as follows. Let $g_t$ be a path of diffeomorphisms of $(\mathbb{C}^n,0)$ such that $g_0$ is the identity.  Then the tangent space $T_fL$ consists of elements of the form
\begin{equation}
\left. \frac{\partial}{\partial t}(f\circ g_t)|_{t=0} = \displaystyle\sum_{i=1}^n \frac{\partial f}{\partial z_i}\cdot\frac{\partial g_i}{\partial t}\right| _{t=0}.
\end{equation}
\noindent In other words, the partial derivatives of $f$ form an $\mathcal{O}_n$-linear basis for $T_fL$, motivating the following important invariants.

Let $I_{\nabla f}\subset \mathcal{O}_n$ be the gradient ideal, generated by the partial derivatives of $f$. Then we define the \emph{local algebra} $A_f := \mathcal{O}_n/I_{\nabla f}$ and its \emph{multiplicity} or \emph{Milnor number} $\mu := \dim A_f$, which are both invariants of the singularity $L$.

Then if $\{\varphi_k\}$ is a monomial basis for $A_f$, we can define a miniversal deformation:
\begin{equation}
F(\lambda) = f + \displaystyle\sum_{k=1}^{\mu} \lambda_k\varphi_k.
\end{equation}
Indeed, the graph of this deformation is a linear subspace of $\mathcal{O}_n$ which is centered at the germ $f\in L$ and is transversal to its orbit.  In particular, this subspace will necessarily intersect every class adjacent to $L$.  To determine these adjacent classes, we restrict to arbitrarily small $\epsilon\in\Lambda$ and use Arnold's algorithm \cite{A1} to classify the possible $F(\epsilon)$.

\subsection{Nonsingular fibers and monodromy}

Let $f:\mathbb{C}^n\rightarrow \mathbb{C}$ be a germ with (isolated) critical point at $0\in\mathbb{C}^n$  of multiplicity $\mu$ and critical value $f(0)=0$.  Let $U$ be a small ball about $0\in\mathbb{C}^n$ and $B$ be a small ball about $0\in\mathbb{C}$.  If the radii of these balls are sufficiently small, the following holds \cite{Mi}.

\begin{thm} \label{Milnor}
For $b\in B' := B\backslash\{0\}$, the level set $X_b = f^{-1}(b)\cap U$ is a nonsingular hypersurface, homotopy equivalent to $\vee^\mu S^{n-1}$.  The level set $X_0 = f^{-1}(0)\cap U$ is nonsingular away from $0$.
\end{thm}

Then $f:X'\rightarrow B'$ (where $X' := f^{-1}(B')\hspace{2pt}\cap\hspace{2pt} U$) is a locally trivial fibration with fiber $X_b \simeq \vee^\mu S^{n-1}$.  Suppose $b_0\in \partial B$ is a noncritical value of $f$, and let $[\gamma] \in \pi_1(B',b_0) \cong \mathbb{Z}$.  Then $\gamma(t)$ lifts to a continuous family of maps $h_t: X_{b_0}\rightarrow X_t$ which can be chosen so that $h_0$ is the identity on $X_{b_0}$ and $h = h_1$ is the identity on $\partial X_{b_0} = f^{-1}(b_0)\cap \partial U$.

The map $h: X_{b_0}\rightarrow X_{b_0}$ is the \emph{monodromy} of $\gamma$.  The induced map on homology,
\begin{equation}
h_\ast : H_{n-1}(X_{b_0})\rightarrow H_{n-1}(X_{b_0}),
\end{equation} 
is the corresponding \emph{monodromy operator}, which is well-defined on the class $[\gamma]$.  If, in addition, $[\gamma]\in \pi_1(B',b_0)$ is a counterclockwise generator, $h_\ast$ is called the \emph{classical-monodromy operator}.

\subsubsection{\sf Vanishing cohomology}
Observe that the (reduced) integral [co]homology is nonzero only in dimension $n-1$, where $H_{n-1}(X_b)\cong \mathbb{Z}^\mu$.  We construct a distinguished basis for this homology group by first considering the simple case where $f$ has a nondegenerate critical point of multiplicity $\mu=1$.

The Morse lemma tells us that in some neighborhood of $0\in\mathbb{C}^n$, there is a coordinate system in which $f(\vec{z}) = z_1^2 + \cdots + z_n^2$.  In this coordinate system, let $S^{n-1} = \{\vec{z} : \|\vec{z}\|^2 = 1, \text{Im}(z_i) = 0\}$ and let $\varphi: [0,1]\rightarrow B$ be a path with $\varphi(0) = b_0$ and $\varphi(1) = 0$.  Then the family of spheres,
\begin{equation}\label{nondeg}
S_t = \sqrt{\varphi(t)}S^{n-1} \subset X_{\varphi(t)},
\end{equation}
\noindent depends continuously on the parameter $t$ and vanishes to the singular point $S_1 = 0\in X_0$.  The sphere $S_0 = \sqrt{b_0}S^{n-1}$ corresponds to a homology class $\Delta \in H_{n-1}(X_{b_0})$, called a \emph{vanishing cycle}.

In the more general case that $f$ has a degenerate critical point of arbitrary multiplicity $\mu$, one can slightly perturb $f$ into a function $f_\epsilon = f + \epsilon g$ with $\mu$ nondegenerate critical points in a small neighborhood of $0\in\mathbb{C}^n$, having distinct critical values $a_i$.  Now consider a system of paths $\varphi_1, \ldots, \varphi_\mu$ with $\varphi_i(0) = b_0$ and $\varphi_i(1) = a_i$.  Suppose that these paths satisfy the following conditions:
\begin{enumerate}
\item The loops formed by traversing $\varphi_i$, followed by a small counterclockwise loop around $a_i$, followed by $\varphi_i^{-1}$ generate $\pi_1(B', b_0)$
\item The paths $\varphi_i$ do not intersect themselves and intersect each other only at $b_0$ for $t=0$
\item The paths are indexed clockwise in $\text{arg}\varphi_i(\epsilon)$
\end{enumerate}
Then, as above, each path $\varphi_i$ determines a distinct vanishing cycle $\Delta_i \in H_{n-1}(X_{b_0})$, and the set $\{\Delta_1,\ldots,\Delta_\mu\}$ form a \emph{distinguished basis of vanishing cycles} for the homology $H_{n-1}(X_{b_0})\cong\mathbb{Z}^\mu$.

\subsection{Mixed Hodge structure in the vanishing cohomology}

For $f:X'\rightarrow B'$ as in Theorem \ref{Milnor}, the $\mu$-dimensional complex vector bundle  $\pi_f^\ast: \mathcal{H}_f^\ast\rightarrow B'$, whose fibers are the complex [co]homology groups $H^{n-1}(X_b;\mathbb{C})$, is called the \emph{vanishing [co]homology bundle} of the singularity $f$.  There is a natural connection $\nabla$ in the vanishing [co]homology bundle, called the \emph{Gauss-Manin connection}, which is defined by covariant derivation $\nabla_b$ along the holomorphic vector field $\frac{\partial}{\partial b}$ on the base $B'$.

We would like to define a mixed Hodge structure in the vanishing cohomology bundle and so review the relevant definitions.  Suppose we have an integer lattice $H_{\mathbb{Z}}$ in a real vector space $H_{\mathbb{R}} = H_{\mathbb{Z}}\otimes_{\mathbb{Z}} \mathbb{R}$.  Let $H = H_{\mathbb{Z}}\otimes_{\mathbb{Z}}\mathbb{C}$ be its complexification.  Then for $k\in\mathbb{Z}$, a \emph{pure Hodge structure of weight $k$} on $H$ is a decomposition: 
\begin{equation}
H = \displaystyle\bigoplus_{p+q=k}H^{p,q},
\end{equation}
into complex subspaces satisfying $H^{p,q}=\overline{H^{q,p}}$, where the bar denotes complex conjugation in $\mathbb{C}$.

Equivalently, we may specify a Hodge structure by a \emph{Hodge filtration}: a finite, decreasing filtration $F^p$ on $H$ satisfying $F^p \oplus \overline{F^{p+1}} = H$.  Indeed, from a Hodge filtration, one can recover a Hodge structure by $H^{p,q} = F^p \cap \overline{F^q}$, and from a Hodge structure, one can recover a Hodge filtration by $F^p = \displaystyle\oplus_{i\geq p}H^{i,k-i}$.  We generalize these notions to a \emph{mixed Hodge structure} on $H$, specified by
\begin{enumerate}
\item A \emph{weight filtration}: a finite, increasing filtration $W_k$ on $H$ which is the complexification of an increasing filtration on $H_\mathbb{Z}\otimes_\mathbb{Z}\mathbb{Q}$,
\item A \emph{Hodge filtration}: a finite, decreasing filtration $F^p$ on $H$,
\end{enumerate}
such that for each $k$, the filtration,
\begin{equation}
F^pgr_k^WH := (F^p\cap W_k + W_{k-1})/W_{k-1},
\end{equation}
satisfies $F^pgr_k^WH \oplus \overline{F^{k-p+1}gr_k^WH} = gr_k^WH$.  That is, $F^pgr_k^WH$ induces a pure Hodge structure of weight $k$ on $gr_k^WH := W_k/W_{k-1}$.

The vanishing cohomology is obtained as the complexification of the integral cohomology of the nonsingular fibers $X_b$.  So to define a mixed Hodge structure in the vanishing cohomology, it remains to specify the relevant weight and Hodge filtrations there.  We follow the construction of \cite{V1, V2, V3}.\\

\subsubsection{\sf Hodge filtration}

\noindent To obtain a Hodge filtration, first consider a holomorphic $(n-1)$-form $\omega$ defined in a neighborhood of $0\in\mathbb{C}^n$.  Since $X_b$ is a complex $(n-1)$-manifold, the restriction $\omega_b = \omega|_{X_b}$ represents a cohomology class $[\omega_b]\in H^{n-1}(X_b, \mathbb{C})$ for all $b\in B'$.  That is, $\omega$ defines a global section $s_\omega: B'\rightarrow \mathcal{H}_f^\ast$, $b\mapsto [\omega_b]$ of the vanishing cohomology bundle.

In the neighborhood of every nonsingular manifold $X_b$, there exists a holomorphic $(n-1)$-form $\omega /df$, with the property that $\omega = df \wedge \omega /df$ in that neighborhood.  As above, the restriction $\omega/df_b = \omega/df|_{X_b}$, called the \emph{residue form}, represents a cohomology class $[\omega/df_b]\in H^{n-1}(X_b, \mathbb{C})$ and defines a global section $\sigma_\omega: B'\rightarrow \mathcal{H}_f^\ast$, $b\mapsto [\omega/df_b]$ of the vanishing cohomology bundle.

The section $\sigma_\omega$ is called a \emph{geometric section}.  For a set of $\mu$ forms that do not satisfy a complex analytic relation, the set of their geometric sections trivalizes the vanishing cohomology bundle, i.e., the corresponding residue forms are a basis in each fiber.

The above sections are holomorphic, meaning that if $\delta_b$ is a cycle in the (integer) homology of the fiber which depends continuously on $b$, (i.e., is covariantly constant via the Gauss-Manin connection), then the map $\sigma_\omega\delta: b\mapsto\int_{\delta_b}\sigma_\omega(b)$ is a holomorphic map $B'\rightarrow\mathbb{C}$.  We consider an asymptotic expansion of such a map around zero.

For example, in the simple case \eqref{nondeg} of a nondegenerate critical point, one can easily see that
\begin{equation}
s_\omega S(b) = \int_{S_b}s_\omega(b) = \int_{B_b} ds_\omega(b) = cb^{n/2} + \cdots.
\end{equation}
$S_b$ is the vanishing sphere, $B_b$ is its interior, and the expansion is proportional to $\text{vol}B_b$ and $d\omega|_0$.

For general $f$ with a (possibly degenerate) critical point at $0$, we can take a set of forms $\omega_1,\ldots,\omega_\mu$ such that their geometric sections trivialize the vanishing cohomology bundle.  Then analysis of the Picard-Fuchs equations of these geometric sections yields the following theorem.

\begin{thm}
Let $\delta_b$ be a continuous family of vanishing cycles over the sector $\theta_0 < \text{arg}b < \theta_1$ in $B'$.  Let $\sigma_\omega$ be a section of the vanishing cohomology.   Then the corresponding integral admits an asymptotic expansion:
\begin{equation}\label{asymp}
\sigma_\omega\delta(b) = \int_{\delta_b}\omega/df_b = \sum_{k,\alpha} \frac{T_{k,\alpha}b^\alpha(\log b)^k}{k!},
\end{equation}
which converges for $b$ sufficiently close to 0.  The numbers $e^{2\pi i \alpha}$ are the eigenvalues of the classical monodromy operator.
\end{thm}

If we fix $\omega$, the coefficients $T_{k,\alpha}$ do not depend on $b$, but they do depend linearly on the section $\delta$, and so determine sections $\tau^\omega_{k,\alpha}$ of the vanishing cohomology bundle via the pairing $\langle \tau^\omega_{k,\alpha}, \delta\rangle = T_{k,\alpha}$ between homology and cohomology.  Thus, we can rewrite the asymptotic expansion \eqref{asymp} as a series expansion of the geometric section:
\begin{equation}
\sigma_\omega = \displaystyle\sum_{k,\alpha} \frac{\tau^\omega_{k,\alpha}b^\alpha(\text{ln} b)^k}{k!}
\end{equation}

This expansion induces a filtration of the vanishing cohomology as follows.  Define
\begin{equation}
\alpha(\omega) := \min\{\alpha : \exists k\geq 0 \text{ such that } \tau^\omega_{k,\alpha}\neq 0\}
\end{equation}
Given a geometric section $\sigma_\omega$, the number $\alpha(\omega)$ is its \emph{order}, and the corresponding expansion:
\begin{equation}
\Sigma_\omega := \displaystyle\sum_k \frac{\tau^\omega_{k,\alpha(\omega)}b^{\alpha(\omega)}(\text{ln} b)^k}{k!}
\end{equation}
is its \emph{principal part}.  Now define a finite, decreasing filtration of $F^p_b$ of the fiber $H^{n-1}(X_b;\mathbb{C})$ by
\begin{equation}
F^p_b := \langle \Sigma_{\omega ,b} : \alpha(\omega) \leq n-p-1\rangle\subseteq H^{n-1}(X_b;\mathbb{C}),
\end{equation}
and the \emph{asymptotic Hodge filtration} filtration of the vanishing cohomology bundle by:
\begin{equation}
F^p := \displaystyle\bigcup_b F^p_b.
\end{equation}

\subsubsection{\sf Weight filtration}

\noindent Suppose we have a nilpotent operator $N$ acting on a finite-dimensional vector space $H$.  Then there is exactly one finite, increasing filtration $W_k$ on $H$ which satisfies:

\begin{enumerate}
\item $N(W_k)\subset W_{k-2}$
\item $N^k: W_{r+k}/W_{r+k-1} \cong W_{r-k}/W_{r-k-1}$ for all $k$
\end{enumerate}

\noindent called the \emph{weight filtration of index} $r$ of $N$.

We obtain a weight filtration in the vanishing cohomology bundle using this construction and the classical-monodromy operator $M$.  As is true for any invertible linear operator, $M$ has a Jordan-Chevalley decomposition $M=M_uM_s$ into commuting unipotent and semisimple parts.  Define a nilpotent operator $N$ to be the logarithm of the unipotent part:
\begin{equation}
N = \displaystyle\sum_i \frac{(-1)^{i+1}(M_u-I)^i}{i}.
\end{equation}
Now for each eigenvalue $\lambda$ of the monodromy operator on $H^{n-1}(X_b;\mathbb{C})$, let $H_{\lambda,b}$ be the corresponding root subspace.  Define a filtration $W_{k,b,\lambda}$ according to the following rules:
\begin{enumerate}
\item If $\lambda = 1$, let $W_{k,b,\lambda}$ be the weight filtration of index $n$ of $N$ on $H_{\lambda,s}$
\item If $\lambda \neq 1$, let $W_{k,b,\lambda}$ be the weight filtration of index $n-1$ of $N$ on $H_{\lambda,s}$
\end{enumerate}

Now define a filtration $W_{k,b}$ of the fiber $H^{n-1}(X_b;\mathbb{C})$ by:
\begin{equation}
W_{k,b} := \displaystyle\bigoplus_\lambda W_{k,b,\lambda}
\end{equation}
and a filtration $W_k$ of the vanishing cohomology bundle by:
\begin{equation}
W_k := \displaystyle\bigcup_b W_{k,b}
\end{equation}
The subbundle $W_k$ is the weight filtration in the vanishing cohomology bundle.  Now we may state the following theorem from \cite{V3}.

\begin{thm}\label{mixedH}
For all $k$ and $p$, the filtrations $W_k$ and $F^p$ are analytic subbundles of the vanishing cohomology bundle, which are invariant under the action of the semisimple part of the monodromy operator.  Furthermore, they specify a mixed Hodge structure in the vanishing cohomology bundle:
\begin{equation}
gr_k^WH = \displaystyle\bigoplus_{p+q=k}H^{p,q},
\end{equation}
\noindent where $H^{p,q} := F^p \cap W_k/F^{p+1}\cap W_k + W_{k-1}$.
\end{thm}

\subsection{Spectrum of a singularity}

In light of Theorem \ref{mixedH}, we are now in a position to define the spectrum of a singularity $f\in K$.\\

Let $f\in K$ be a singularity.  If $\lambda$ is an eigenvalue of the semisimple part of the classical-monodromy operator on $H^{p,q}$, one can associate to $f$ the set of $\mu$ rational numbers:
\begin{equation}
\{n-1-l_p\lambda\}\text{, }\hspace{10pt}\text{where }\hspace{10pt}l_p\lambda := \log(\lambda/2\pi i)\hspace{10pt}\text{ and }\hspace{10pt}\text{Re}(l_p\lambda) = p.
\end{equation}
This (unordered) set of numbers is the \emph{spectrum} of the singularity $K$.

To see what the spectrum of a singularity $f$ has to do with adjacencies to $f$, we first construct a fibration, analagous to the fibration $f:X'\rightarrow B'$.  Choose a miniversal deformation:
\begin{equation}
F(z,\lambda) = f(z) + \displaystyle\sum_{i=0}^{\mu-1}\lambda_i\varphi_i(z)
\end{equation}
\noindent where $\lambda\in \mathbb{C}^\mu$ and $\varphi_0 = 1$.  As before, we consider sufficiently small ball $U$ about $0\in\mathbb{C}^n$ and a another small ball, this time $\Lambda$ about $0\in\mathbb{C}^\mu$.

For $\lambda\in\Lambda$, define the level set $V_\lambda := \{z\in U : F(z,\lambda)=0\}$ and the hypersurface $V := \{(z,\lambda)\in U\times\Lambda : F(z,\lambda)=0\}$. Let $\Sigma\subset\Lambda$ be the set of values of $\lambda$ for which $V_\lambda$ is singular, called the \emph{level bifurcation set}.  Let $\pi_\Lambda: V\rightarrow \Lambda$ be the restriction of the canonical projection, called the \emph{Whitney map}.  Finally, let $\Lambda' := \Lambda\backslash\Sigma$ and $V' := \pi^{-1}_\Lambda(\Lambda')$.  The locally trivial fibration $\pi_\Lambda: V'\rightarrow \Lambda'$ with fiber $V_\lambda$ over $\lambda\in\Lambda'$ is the \emph{Milnor fibration} of $f$.

Observe that the fibration $f:X'\rightarrow B'$ can be embedded in the Milnor fibration by identifying $B'$ with the $\lambda_0$-axis in the base $\Lambda'$ (recall that $\varphi_0 = 1$).  Furthermore, we can repeat the constructions outlined above for the Milnor fibration and then ask how the spectrum varies as we vary the deformation parameter $\lambda$ in an infinitesimal neighborhood of 0.  This leads to observations on the \emph{semicontinuity} of the spectrum, including the following \cite{A2}.

\begin{thm}\label{adjspec}
Suppose that a critical point of type $L$ has (ordered) spectrum $\alpha_1\leq\cdots\leq\alpha_\mu$ and a critical point of type $L'$ has spectrum $\alpha'_1\leq\cdots\leq\alpha'_{\mu'}$ where $\mu' \leq \mu$.  Then a necessary condition for the adjacency $L\rightarrow L'$ is that the spectra be adjacent in the sense that $\alpha_i \leq \alpha'_i$.
\end{thm}

\newpage
\appendix
\setcounter{equation}{0}
\section{DAHA-Jones formulas}
\label{sec:AandD}

\noindent {\bf Type $A$.} The formulas for $\widetilde{J\!D}^{A_n}(b)$, can be readily obtained
from the following well-known type-$A$ super-polynomials
$H\!D_{r,s}^{A} (b;q,t,a)$ upon the substitution $a=-t^{n+1}$.
We will need only $A_6$ here, which corresponds to $a=-t^7$:
\begin{gather}\label{Asup}
H\!D_{3,2}^{A}(\omega_1)=1 + a q + q t, \\
H\!D_{5,2}^{A}(\omega_1)=1+ q t+q^2 t^2+a \bigl(q+q^2 t\bigr),\\
H\!D_{4,3}^{A}(\omega_1)=1+ a^2 q^3+q t+q^2 t+q^2 t^2+q^3 t^3 + a \bigl(q+q^2+q^2 t+q^3 t
+q^3 t^2\bigr).
\end{gather}

The simplest colored formulas for the super-polynomials of type $A$, defined for $\omega_2$, are known from \cite{GS,FGS} and \cite{ChJ}.  They play an important role for the super-polynomials of the pair $(E_6,\omega_1)$, in spite of the fact that this weight is non-minuscule.
\begin{gather}
H\!D^{A}_{3,2}(\omega_2;q,t,a)\ =\\1+\frac{a^2 q^2}{t}+q t+q t^2+q^2 t^4+a \bigl(q+\frac{q}{t}
+q^2 t+q^2 t^2\bigr),\notag
\end{gather}

\begin{equation}
H\!D^{A}_{5,2}(\omega_2;q,t,a)\ =
\end{equation}
\renewcommand{\baselinestretch}{0.5}
{\normalsize
\noindent\(
1+q t+q t^2+q^2 t^2+q^2 t^3+q^2 t^4+q^3 t^5+q^3 t^6+q^4 t^8
+a^2 \bigl(q^3+\frac{q^2}{t}+q^3 t+q^4 t^3\bigr)
+a \bigl(q+q^2+\frac{q}{t}+2 q^2 t
+q^2 t^2+q^3 t^2+2 q^3 t^3+q^3 t^4+q^4 t^5+q^4 t^6\bigr),
\)}
\renewcommand{\baselinestretch}{1.2}

\begin{equation}
H\!D^{A}_{4,3}(\omega_2;q,t,a)\ =
\end{equation}
\renewcommand{\baselinestretch}{0.5}
{\small
\noindent\(
1+\frac{a^4 q^6}{t^2}+q t+q^2 t+q t^2+2 q^2 t^2
+q^2 t^3+2 q^3 t^3+q^2 t^4+2 q^3 t^4+q^4 t^4+q^3 t^5
+q^4 t^5+q^3 t^6+2 q^4 t^6+q^4 t^7+q^5 t^7+q^4 t^8
+q^5 t^8+q^5 t^9+q^5 t^{10}+q^6 t^{12}
+a^3 \bigl(q^5+q^6+\frac{q^4}{t^2}+\frac{q^5}{t^2}
+\frac{q^4}{t}+\frac{q^5}{t}+q^5 t+q^6 t+q^6 t^2
+q^6 t^3\bigr)+a^2 \bigl(2 q^3+2 q^4+q^5
+\frac{q^3}{t^2}+\frac{q^2}{t}+\frac{2 q^3}{t}
+\frac{q^4}{t}+q^3 t+4 q^4 t+2 q^5 t+2 q^4 t^2
+3 q^5 t^2+q^4 t^3+3 q^5 t^3+q^6 t^3+2 q^5 t^4+q^6 t^4
+q^5 t^5+2 q^6 t^5+q^6 t^6+q^6 t^7\bigr)
+a \bigl(q+2 q^2+q^3+\frac{q}{t}+\frac{q^2}{t}
+2 q^2 t+4 q^3 t+q^4 t+q^2 t^2+4 q^3 t^2+2 q^4 t^2
+2 q^3 t^3+4 q^4 t^3+q^5 t^3+q^3 t^4+4 q^4 t^4+2 q^5 t^4
+2 q^4 t^5+3 q^5 t^5+q^4 t^6+3 q^5 t^6+2 q^5 t^7
+q^6 t^7+q^5 t^8+q^6 t^8+q^6 t^9+q^6 t^{10}\bigr).
\)}
\renewcommand{\baselinestretch}{1.2}
\medskip

More specifically, we will need the values of these super-polynomials at $t=1$:
\begin{gather}
H\!D^{A}_{3,2}(\omega_2;q,t=1,a)\ =\
(1+q+a q)^2,\\
H\!D^{A}_{5,2}(\omega_2;q,t=1,a)\ =\
(1 + q + a q + q^2 + a q^2)^2,\\
H\!D^{A}_{4,3}(\omega_2;q,t=1,a)=
(1 + q + a q + 2 q^2 + 2 a q^2 + q^3 + 2 a q^3 + a^2 q^3)^2.
\end{gather}
For instance, the corresponding {\em dimensions}
$H\!D^{A}(q=1,t=1,a=1)$  are
$9, 25, 121$.
\medskip

\noindent {\bf Type $D$.} We will need the following DAHA-Jones
polynomials of type $D_5$ for $\omega_1$
(which is minuscule):
\begin{gather}
\widetilde{J\!D}^{D_5}_{3,2}(\omega_1;q,t)\ =\\
1+q t+q t^4-q t^5-q t^8+q^2 t^8-q^2 t^9-q^2 t^{12}+q^2 t^{13},\notag
\end{gather}
\begin{equation}
\widetilde{J\!D}^{D_5}_{5,2}(\omega_1;q,t)\ =
\end{equation}
\renewcommand{\baselinestretch}{0.5}
{\normalsize
\noindent\(
1+q t+q^2 t^2+q t^4-q t^5+q^2 t^5-q^2 t^6-q t^8
+q^2 t^8-2 q^2 t^9+q^3 t^9-q^3 t^{10}-q^2 t^{12}
+q^3 t^{12}+q^2 t^{13}-2 q^3 t^{13}+q^3 t^{14}
-q^3 t^{16}+q^4 t^{16}+q^3 t^{17}-q^4 t^{17}-q^4 t^{20}+q^4 t^{21},
\)}
\renewcommand{\baselinestretch}{1.2}
\medskip

\begin{equation}
\widetilde{J\!D}^{D_5}_{4,3}(\omega_1;q,t)\ =
\end{equation}
\renewcommand{\baselinestretch}{0.5}
{\normalsize
\noindent\(
1+q t+q^2 t^2+q t^4-q t^5+q^2 t^5-q^2 t^6-q t^8+q^2 t^8-2 q^2 t^9
+q^3 t^9-q^3 t^{10}-q^2 t^{12}+q^3 t^{12}+q^2 t^{13}-2 q^3 t^{13}
+q^3 t^{14}-q^3 t^{16}+q^4 t^{16}
+q^3 t^{17}-q^4 t^{17}-q^4 t^{20}+q^4 t^{21}.
\)}
\renewcommand{\baselinestretch}{1.2}
\medskip

We will also need the super-polynomials for
the case when the {\em last} fundamental weight is
taken for $D_n (n\ge 4)$:
\begin{gather}\label{Dsuplast}
\widehat{H\!D}_{3,2}^{D}(\omega_n)=1 + a qt^6 + q t^3,\\
\widehat{H\!D}_{5,2}^{D}(\omega_n)=1+ q t^3+ q^2 t^6+a \bigl(q t^6+q^2 t^9\bigr),\\
\widehat{H\!D}_{4,3}^{D}(\omega_n)=1+  a^2 q^3 t^{14} +q t^3+q^2 t^5+q^2 t^6+a \bigl(q t^6+q^2+q^2 t^8+q^2 t^9 +q^3 t^{11}
+q^3 t^{12}\bigr),
\end{gather}
where the relevant specializations are
\begin{equation}
\widehat{H\!D}^{D}(q, t, a\mapsto -t^{n-4}) = 
\widetilde{J\!D}^{D_n}(\omega_n\,;\, q,t).
\end{equation}
The DAHA-superpolynomials and DAHA-Jones polynomials for $\omega_{n-1}$ are identical to those for $\omega_n$.

Interestingly, these super-polynomials are
related to those for ($A$, $\omega_1$):
\begin{align}\label{hatda}
\widehat{H\!D}_{r,s}^{D}(\omega_n;q,t,a)=
H\!D^A_{r,s}(\omega_1;q\mapsto t q^2,t, a\mapsto a t^4),
\end{align}
so we have essentially similar ``stable theories" for the pairs
$(A_{n-1}, \omega_1)$ and $(D_n, \omega_n)$.
\medskip

\noindent {\bf Type $E_6$.} We will need the DAHA-Jones polynomials for the minuscule weight $\omega_1$:
\begin{gather}
\widetilde{J\!D}^{E_6}_{3,2}(\omega_1;q,t) =\\
1+q \bigl(t+t^4-t^9-t^{12}\bigr)+q^2 \bigl(t^8-t^{13}-t^{16}
+t^{21}\bigr),\notag
\end{gather}

\begin{equation}
\widetilde{J\!D}^{E_6}_{5,2}(\omega_1;q,t) =
\end{equation}
\renewcommand{\baselinestretch}{0.5}
{\normalsize
\noindent\(
1+q \bigl(t+t^4-t^9-t^{12}\bigr)
+q^2 \bigl(t^2+t^5+t^8-t^{10}-2 t^{13}-t^{16}+t^{21}\bigr)
+q^3 \bigl(t^9+t^{12}-t^{14}-2 t^{17}-t^{20}+t^{22}+t^{25}\bigr)
+q^4 \bigl(t^{16}-t^{21}-t^{24}+t^{29}\bigr),
\)}
\renewcommand{\baselinestretch}{1.2}

\begin{equation}
\widetilde{J\!D}^{E_6}_{4,3}(\omega_1;q,t) =
\end{equation}
\renewcommand{\baselinestretch}{0.5}
{\small
\noindent\(
1+q \bigl(t+t^4-t^9-t^{12}\bigr)+q^2 \bigl(t+t^2+t^4+t^5+t^8-t^9
-t^{10}-t^{12}-2 t^{13}-t^{16}+t^{21}\bigr)
+q^3 \bigl(t^3+t^5+t^6+t^8+t^9-t^{10}-t^{11}+t^{12}
-3 t^{13}-2 t^{14}-2 t^{16}-2 t^{17}+t^{18}-t^{20}+2 t^{21}
+t^{22}+t^{24}+t^{25}\bigr)+q^4 \bigl(t^8+t^9+t^{10}+t^{12}
-t^{14}-t^{15}-3 t^{17}-2 t^{18}-2 t^{20}-t^{21}+2 t^{22}
+t^{23}-t^{24}+3 t^{25}+t^{26}+t^{28}+t^{29}-t^{30}-t^{33}\bigr)
+q^5 \bigl(t^{13}+t^{16}-t^{18}-3 t^{21}-2 t^{24}+2 t^{26}
+3 t^{29}+t^{32}-t^{34}-t^{37}\bigr)+q^6 \bigl(t^{24}-t^{25}
-t^{28}+t^{30}-t^{32}
+2 t^{33}-t^{34}+t^{36}-t^{38}-t^{41}+t^{42}\bigr).
\)}
\renewcommand{\baselinestretch}{1.2}
\medskip

\noindent The next series of DAHA-Jones polynomials will be for $\omega_6$ (minuscule):
\begin{gather}
\widetilde{J\!D}^{E_6}_{3,2}(\omega_6;q,t) =\\
1 + q (t^{5} + t^{8} - t^{9} - t^{12}) + q^2 (t^{16} - t^{17}
- t^{20} + t^{21}),\notag
\end{gather}

\begin{equation}
\widetilde{J\!D}^{E_6}_{5,2}(\omega_6;q,t) =
\end{equation}
\renewcommand{\baselinestretch}{0.5}
{\normalsize
\noindent\(
1+q \bigl(t^5+t^8-t^9-t^{12}\bigr)+q^2 \bigl(t^{10}+t^{13}-t^{14}
+t^{16}-2 t^{17}-t^{20}+t^{21}\bigr)
+q^3 \bigl(t^{21}-t^{22}+t^{24}-2 t^{25}
+t^{26}-t^{28}+t^{29}\bigr)
+q^4 \bigl(t^{32}-t^{33}-t^{36}+t^{37}\bigr),
\)}
\renewcommand{\baselinestretch}{1.2}

\begin{equation}
\widetilde{J\!D}^{E_6}_{4,3}(\omega_6;q,t) =
\end{equation}
\renewcommand{\baselinestretch}{0.5}
{\normalsize
\noindent\(
1+q \bigl(t^5+t^8-t^9-t^{12}\bigr)
+q^2 \bigl(t^9+t^{10}+t^{12}-t^{14}-2 t^{17}-t^{20}+t^{21}\bigr)
+q^3 \bigl(t^{15}+t^{17}-t^{19}+t^{20}
-2 t^{21}-t^{22}-t^{24}+t^{26}+t^{29}\bigr)
+q^4 \bigl(t^{24}-t^{27}-t^{29}+t^{31}-t^{32}+t^{33}\bigr).
\)}
\renewcommand{\baselinestretch}{1.2}
\newpage

\setcounter{equation}{0}
\section{Figures}\label{fig}

\begin{figure}[htbp]
\begin{centering}
\includegraphics[width=6in]{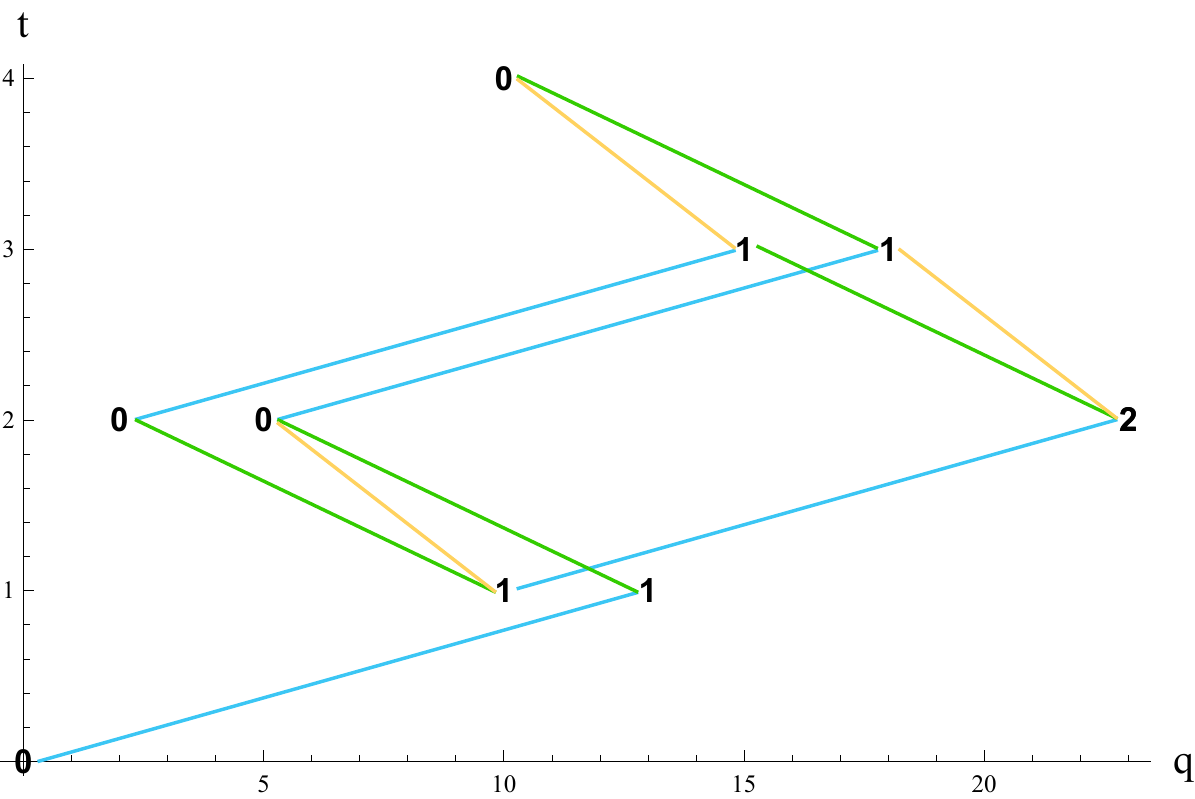}
\caption{Differentials for $T^{3,2}$}
\label{fig:T23ALL}
\end{centering}
\end{figure}

This appendix contains diagrams which depict our proposals for $H^{\frak e_6, \mathbf{27}}$ in Section \ref{sec:hyper}.  We use QG-conventions; see \eqref{CDconv}.  In particular, Figure \ref{fig:T23ALL} corresponds to our proposal for $T^{3,2}$, figure \ref{fig:T25ALL} corresponds to our proposal for $T^{5,2}$, and figure \ref{fig:T25ALL} corresponds to our proposal for $T^{4,3}$.

In each figure, a monomial $q^it^ju^k$ corresponds to the number $k$ placed on the diagram in position $(i,j)$, i.e., with $x$-coordinate $i$ and $y$-coordinate $j$.  The differentials are depicted by line segments connecting pairs of monomials, color-coded as follows. 

\begin{equation}\label{diffcolor}
    \begin{tabular}{| l | c | l |}
    \hline
    $\mathfrak{g},V$ & color & deg$(d_{\mathfrak{g},V})$ \\ \hline \hline
    $\mathfrak{e}_6,\mathbf{27}$ & -- & $(0,-1,1)$  \\ \hline
    $\mathfrak{d}_5, \mathbf{10}$ & {\color{red} Red} & $(4,-1,1)$ \\ \hline
    $\mathfrak{a}_6, \mathbf{7}$ & {\color{YellowOrange} Yellow} & $(5,-1,1)$  \\ \hline
    canceling & {\color{Green} Green} & $(8,-1,1)$  \\ \hline
    canceling & {\color{CornflowerBlue} Blue} & $(13,1,1)$  \\
    \hline
    \end{tabular}
\end{equation}
Observe that while the differential corresponding to $(\mathfrak{d}_5, \mathbf{10})$ only appears in the diagram for $H^{\frak e_6, \mathbf{27}}(T^{4,3})$, that structure still exists as a \emph{specialization} in the other two cases; see Section \ref{sec:hyper}.

\newpage

\begin{centering}
\begin{landscape}
\begin{figure}
\centering
\includegraphics[width=8.5in]{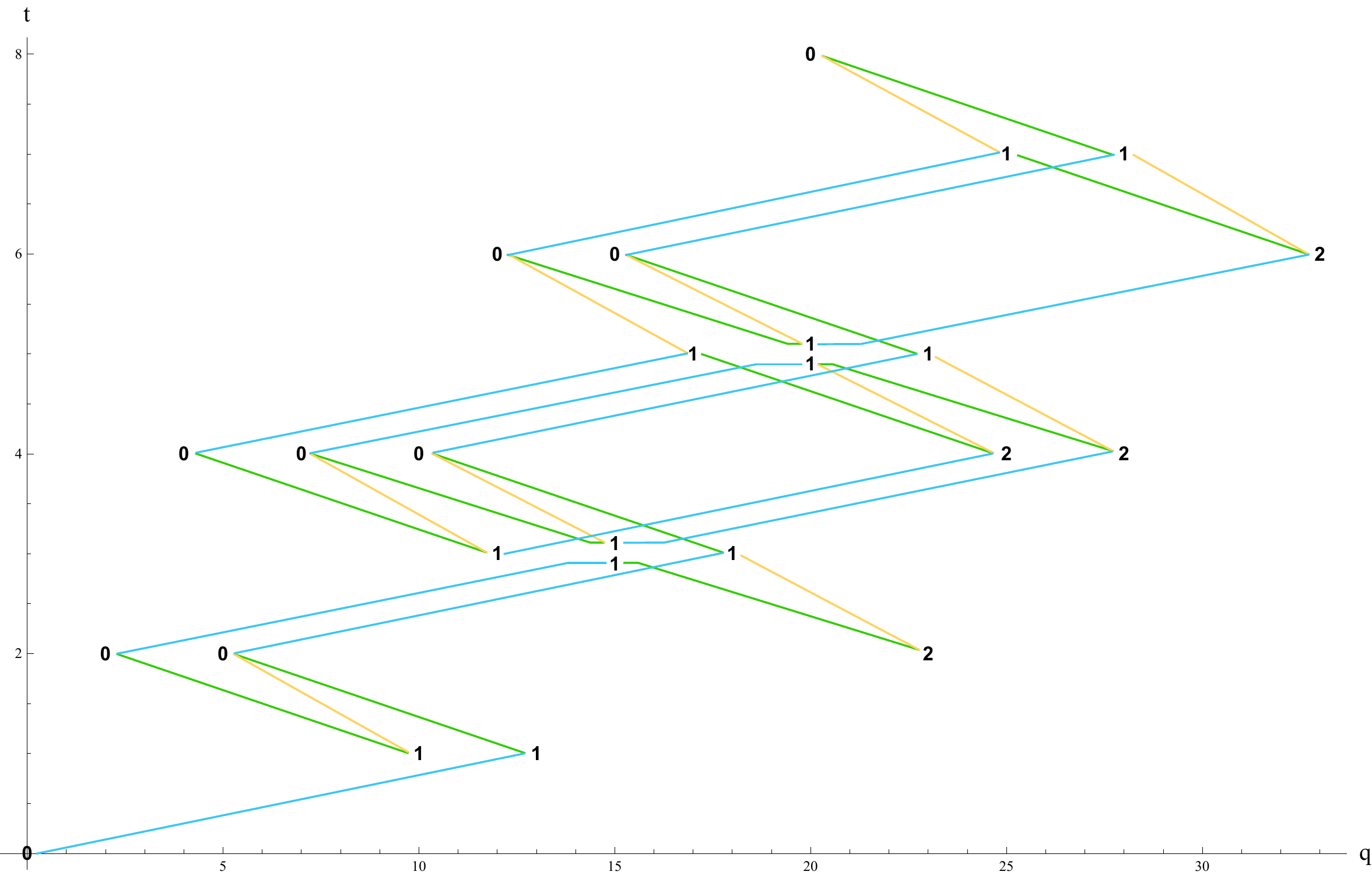}
\caption{Differentials for $T^{5,2}$}
\label{fig:T25ALL}
\end{figure}
\end{landscape}
\end{centering}
\newpage


\begin{landscape}
\begin{center}
\begin{figure}
\centering
\includegraphics[width=8.5in]{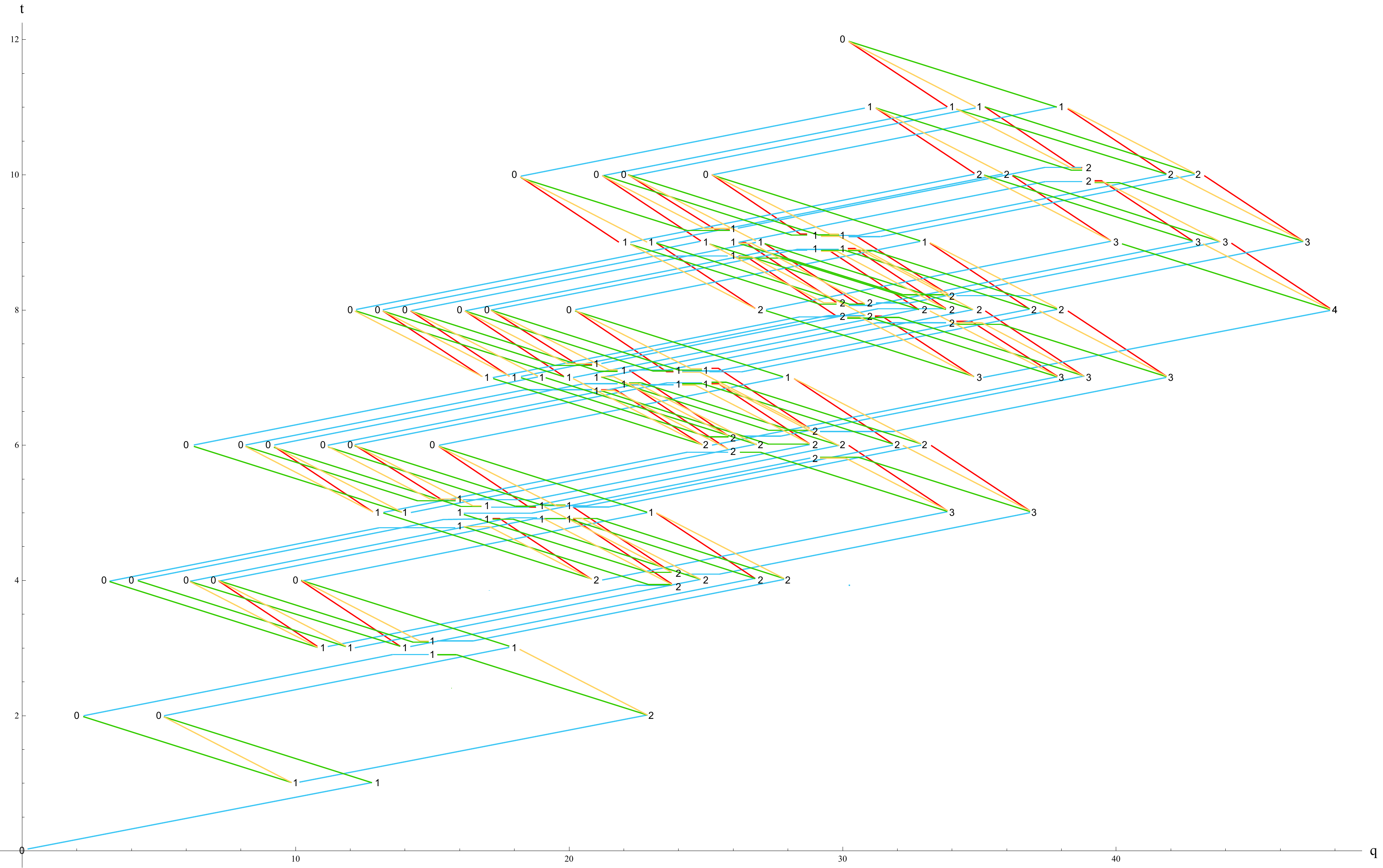}
\caption{Differentials for $T^{4,3}$}
\label{fig:T34ALL}
\end{figure}
\end{center}
\end{landscape}
\newpage

\setcounter{equation}{0}
\section{Adjacencies and spectra}\label{app:Adj}

This appendix contains the adjacency tree to $Z_{3,0}$, computed as outlined in section \ref{sec:Adj}.  This tree displays only those adjacencies (``arrows") that arise in the classification of singularities by their jets \cite{A1}, though there are other internal adjacencies.  Observe that as a direct consequence of the definition of adjacency, this tree is transitive in that $A\rightarrow B\rightarrow C$ implies $A\rightarrow C$.

One can check this list using the adjacency of the spectra, which are also listed.  There are many ways to compute the spectrum of a singularity, and we will outline one method here.  Suppose $f\in\mathcal{O}_n$ with Taylor expansion $f = \sum a_{\mathbf{k}}\mathbf{z^k}$.  Then we can take the set:
\begin{equation}
\text{supp} f = \{\mathbf{k}\in \mathbb{N}^n_{\geq 0} : a_{\mathbf{k}}\neq 0\}
\end{equation}
Now we let define a subset of $\mathbb{R}^n_+$ by:
\begin{equation}
G(f) = \displaystyle\bigcup_{\mathbf{k}\in\text{supp} f} \{\mathbf{k} + \mathbb{R}^n_+\}
\end{equation}
The convex hull of $G(f)$ constitutes the \textit{Newton polyhedron} of $f$, and the union of the compact faces of the Newton polyhedron is the \textit{Newton diagram} $\Gamma(f)$ of $f$.

A Newton diagram induces a decreasing filtration on power series as follows.  If we assume that any monomial contained in the Newton diagram is quasihomogeneous of degree 1, then each face $e_i\in \Gamma(f)$ determines a set of weights $\mathbf{\nu}_i$ such that $\langle \mathbf{j}, \mathbf{\nu}_i\rangle = 1$ for all $\mathbf{z^j}\in e_i$.  We can then define the Newton degree of an arbitrary monomial by:
\begin{equation}
\text{deg}\mathbf{z^k} = \min_i \langle\mathbf{k}, \mathbf{\nu}_i\rangle.
\end{equation}
Then if every monomial in a power series has Newton degree greater than or equal to $d$, that power series belongs to the $d^{th}$ subspace of the Newton filtration.

The Newton filtration also descends to forms, e.g., the Newton order of the form $\mathbf{z^k}dz_1\wedge\cdots\wedge dz_n$ coincides with the Newton order of the monomial $\mathbf{z^k}z_1\cdots z_n$.  Furthermore, the Newton filtration on forms coincides with the Hodge filtration after a shift of indices, and one can show that for an appropriate set of monomials (ones whose corresponding forms trivialize the vanishing cohomology bundle), the spectrum coincides with the set of numbers:
\begin{equation}
\min_i\langle \textbf{k}+\textbf{1}, \boldsymbol\nu_i\rangle - 1,
\end{equation}
for those monomials, which can often be taken to be a basis for the local algebra or, using the symmetry of the spectrum about $\frac{n}{2}-1$, a set of \textit{subdiagrammatic monomials}--those $\mathbf{z^k}$ for which $\mathbf{k} + \mathbf{1}$ does not belong to the interior of the Newton polyhedron.  
\newpage

\noindent The following table lists the singularities adjacent to $Z_{3,0}$ and their normal forms, relevant deformations, and Milnor numbers.

\footnotesize
\bigskip
\begin{center}
    \begin{tabular}{| l | l | l | l |}
    \hline
    \textbf{Singularity} & \textbf{Normal Form}$^1$ & $\mathbf{\Delta W}/\epsilon$ & $\mathbf{\mu}$ \\ \hline \hline
    $Z_{3,0}$ & $x^3y + dx^2y^5 + \mathbf{a}_4xy^{10} + y^{13}$ & $--$ & $27$ \\ \hline
    $A_k$, $1\leq k\leq 12$ & $x^{k+1}$ & $x^2 + y^{k+1}$ & $k$ \\ \hline
    $D_k$, $4\leq k\leq 14$ & $x^2y + y^{k-1}$ & $x^2y + y^{k-1}$ & $k$  \\ \hline
    $E_{6}$ & $x^3+y^4$ & $x^3 + y^4$ & $6$  \\ \hline
    $E_{7}$ & $x^3 + xy^3$ & $x^3 + xy^3$ & $7$ \\ \hline
    $E_{8}$ & $x^3+y^5$ & $x^3 + y^5$ & $8$  \\ \hline
    $J_{2,0}$ & $x^3 + bx^2y^2 + y^6$ & $x^3 + x^2y^2 + y^6$ & $10$  \\ \hline
    $J_{2,p}$, $1\leq p\leq 7$ & $x^3 + x^2y^2 + ay^{6+p}$ & $x^3 + x^2y^2 + y^{6+p}$ & $10+p$  \\ \hline
    $E_{12}$ & $x^3+y^7 + \mathbf{a}_2xy^5$ & $x^3 + y^{7}$ & $12$  \\ \hline
    $E_{13}$ & $x^3 + xy^5 + \mathbf{a}_2y^8$ & $x^3 + xy^5$ & $13$ \\ \hline
    $E_{14}$ & $x^3+ y^8 +\mathbf{a}_2xy^6$ & $x^3 + y^8$ & $14$  \\ \hline
    $J_{3,0}$ & $x^3 + bx^2y^3 + y^9 + axy^7$ & $x^3 + x^2y^3 + y^9$ & $16$  \\ \hline
    $J_{3,p}$, $1\leq p\leq 4$ & $x^3 + x^2y^3 + \mathbf{a}_3y^{9+p}$ & $x^3 + x^2y^3 + y^{9+p}$ & $16+p$  \\ \hline
    $E_{18}$ & $x^3 + y^{10} + \mathbf{a}_3xy^7$ & $x^3 + y^{10}$ & $18$  \\ \hline
    $E_{19}$ & $x^3 + xy^7 + \mathbf{a}_3y^{11}$ & $x^3 + xy^7$ & $19$ \\ \hline
    $E_{20}$ & $x^3+ y^{11} +\mathbf{a}_3xy^8$ & $x^3 + y^{11}$ & $20$ \\ \hline
    $J_{4,0}$ & $x^3 + bx^2y^4 + y^{12} + \mathbf{a}_3xy^9$ & $x^3 + x^2y^4 + y^{12}$ & $22$  \\ \hline
    $J_{4,1}$ & $x^3 + x^2y^4 + \mathbf{a}_4y^{13}$ & $x^3 + x^2y^4$ & $23$  \\ \hline
    $J_{4,2}$ & $x^3 + x^2y^4 + \mathbf{a}_4y^{14}$ & $x^3 + xy^9 + x^2y^4$ & $24$  \\ \hline
    $E_{24}$ & $x^3 + y^{13} + \mathbf{a}_4xy^9$ & $x^3 + y^{10}$ & $24$  \\ \hline
    $E_{25}$ & $x^3 + xy^9 + \mathbf{a}_4y^{14}$ & $x^3 + xy^9$ & $25$ \\ \hline
    $X_{1,0}$ & $x^4 + ax^2y^2 + y^4$, $a \neq 4$ & $y^4 + xy^3 + x^2y^2$ & $9$  \\ \hline
    $X_{1,p}$, $1\leq p\leq 9$ & $x^4 + x^2y^2 + ay^{4+p}$, $a\neq 0$ & $x^2y^2 + y^{4+p}$ & $9+p$  \\ \hline
    $Z_{11}$ & $x^3y + y^5 + axy^4$ & $y^5$ & $11$  \\ \hline
    $Z_{12}$ & $x^3y + xy^4 + ax^2y^3$ & $xy^4$ & $12$ \\ \hline
    $Z_{13}$ & $x^3y + y^6 + axy^5$ & $y^6$ & $13$  \\ \hline
    $Z_{1,0}$ & $x^3y + dx^2y^3 + axy^6 + y^7$ & $x^2y^3 + y^7$ & $15$  \\ \hline
    $Z_{1,p}$, $1\leq p\leq 6$ & $x^3y + x^2y^3 + \mathbf{a}_3y^{7+p}$ & $x^2y^3 + y^{7+p}$ & $15+p$  \\ \hline
    $Z_{17}$ & $x^3y + y^8 + \mathbf{a}_3xy^6$ & $y^8$ & $17$  \\ \hline
    $Z_{18}$ & $x^3y + xy^6 + \mathbf{a}_3y^9$ & $xy^6$ & $18$ \\ \hline
    $Z_{19}$ & $x^3y + y^9 + \mathbf{a}_3xy^7$ & $y^9$ & $19$  \\ \hline
    $Z_{2,0}$ & $x^3y + dx^2y^4 + \mathbf{a}_3xy^8 + y^{10}$ & $x^2y^4 + y^{10}$ & $21$  \\ \hline
    $Z_{2,p}$, $1\leq p\leq 4$ & $x^3y + x^2y^4 + \mathbf{a}_4y^{9+p}$ & $x^2y^4 + y^{9+p}$ & $21+p$  \\ \hline
    $Z_{23}$ & $x^3y + y^{11} + \mathbf{a}_4xy^8$ & $y^{11}$ & $23$  \\ \hline
    $Z_{24}$ & $x^3y + xy^8 + \mathbf{a}_4y^{12}$ & $xy^8$ & $24$ \\ \hline
    $Z_{25}$ & $x^3y + y^{12} + \mathbf{a}_4xy^9$ & $y^{12}$ & $25$  \\
    \hline
    \end{tabular}
\end{center}
\let\thefootnote\relax\footnote{$^1$Here we have that $\mathbf{a}_k := a_0 + \cdots + a_{k-2}y^{k-2}, \mathbf{a}_1 := 0$.}

\newpage
\normalsize
\bigskip
\noindent The following table lists the spectra of the singularities which are adjacent to $Z_{3,0}$.  Observe that, by Theorem \ref{adjspec}, it supports our list of adjacencies.
\bigskip
\small
\begin{center}
    \begin{tabular}{| l | p{13.5cm} |}
    \hline
     & \textbf{Spectrum}  \\ \hline \hline
    $Z_{3,0}$ & $\scriptscriptstyle{\frac{-3}{26} , \frac{-1}{26} , \frac{1}{26} , \frac{3}{26} , \frac{5}{26} , \frac{5}{26} , \frac{7}{26} , \frac{7}{26} , \frac{9}{26} , \frac{9}{26} , \frac{11}{26} , \frac{11}{26} , \frac{13}{26} , \mathbf{\frac{13}{26}}, \frac{13}{26} , \frac{15}{26} , \frac{15}{26} , \frac{17}{26} , \frac{17}{26} , \frac{19}{26} , \frac{19}{26} , \frac{21}{26} , \frac{21}{26} , \frac{23}{26} , \frac{25}{26} , \frac{27}{26} , \frac{29}{26}}$ \\ \hline
    $A_k$ & $\scriptscriptstyle{\frac{2}{2k+2} , \frac{4}{2k+2} , \cdots , \frac{2k}{2k+2}}$ \\ \hline
    $D_k$ & $\scriptscriptstyle{\left\{\frac{1}{2k-2} , \frac{3}{2k-2} , \cdots , \frac{2k-3}{2k-2}\right\}}$ $\cup$ $\scriptscriptstyle{\left\{\frac{k-1}{2k-2}\right\}}$  \\ \hline
    $E_{6}$ & $\scriptscriptstyle{\frac{1}{12} , \frac{4}{12} , \frac{5}{12} , \frac{7}{12} , \frac{8}{12} , \frac{11}{12}}$ \\ \hline
    $E_{7}$ & $\scriptscriptstyle{\frac{1}{18} , \frac{5}{18} , \frac{7}{18} , \mathbf{\frac{9}{18}} , \frac{11}{18} , \frac{13}{18} , \frac{17}{18}}$ \\ \hline
    $E_{8}$ & $\scriptscriptstyle{\frac{1}{30} , \frac{7}{30} , \frac{11}{30} , \frac{13}{30} , \frac{17}{30} , \frac{19}{30} , \frac{23}{30} , \frac{29}{30}}$ \\ \hline
    $J_{2,p}$ & $\scriptscriptstyle{\left\{\frac{0}{6(p+6)} , \frac{6}{6(p+6)} , \cdots , \frac{6(p+6)}{6(p+6)}\right\}}$ $\cup$ $\scriptscriptstyle{\left\{\frac{2(p+6)}{6(p+6)} , \frac{3(p+6)}{6(p+6)} , \frac{4(p+6)}{6(p+6)}\right\}}$ \\ \hline
    $E_{12}$ & $\scriptscriptstyle{\frac{-1}{42} , \frac{5}{42} , \frac{11}{42} , \frac{13}{42} , \frac{17}{42} , \frac{19}{42} , \frac{23}{42} , \frac{25}{42} , \frac{29}{42} , \frac{31}{42} , \frac{37}{42} , \frac{43}{42}}$ \\ \hline
    $E_{13}$ & $\scriptscriptstyle{\frac{-1}{30} , \frac{3}{30} , \frac{7}{30} , \frac{9}{30} , \frac{11}{30} , \frac{13}{30} , \mathbf{\frac{15}{30}} , \frac{17}{30} , \frac{19}{30} , \frac{21}{30} , \frac{23}{30} , \frac{27}{30} , \frac{31}{30}}$ \\ \hline
    $E_{14}$ & $\scriptscriptstyle{\frac{-1}{24} , \frac{2}{24} , \frac{5}{24} , \frac{7}{24} , \frac{8}{24} , \frac{10}{24} , \frac{11}{24} , \frac{13}{24} , \frac{14}{24} , \frac{16}{24} , \frac{17}{24} , \frac{19}{24} , \frac{22}{24} , \frac{25}{24}}$  \\ \hline
    $J_{3,p}$ & $\scriptscriptstyle{\left\{\frac{9}{18(p+9)} , \frac{27}{18(p+9)} ,\cdots , \frac{9(2p+17)}{18(p+9)}\right\}}$ $\cup$ $\scriptscriptstyle{\left\{\frac{-(p+9)}{18(p+9)} , \frac{5(p+9)}{18(p+9)}, \frac{7(p+9)}{18(p+9)},\frac{9(p+9)}{18(p+9)} , \frac{11(p+9)}{18(p+9)}, \frac{13(p+9)}{18(p+9)},\frac{19(p+9)}{18(p+9)}\right\}}$ \\ \hline
    $E_{18}$ & $\scriptscriptstyle{\frac{-2}{30} , \frac{1}{30} , \frac{4}{30} , \frac{7}{30} , \frac{8}{30} , \frac{10}{30} , \frac{11}{30} , \frac{13}{30} , \frac{14}{30} , \frac{16}{30} , \frac{17}{30} , \frac{19}{30} , \frac{20}{30} , \frac{22}{30} , \frac{23}{30} , \frac{26}{30} , \frac{29}{30} , \frac{32}{30}}$ \\ \hline
    $E_{19}$ & $\scriptscriptstyle{\frac{-3}{42} , \frac{1}{42} , \frac{5}{42} , \frac{9}{42} , \frac{11}{42} , \frac{13}{42} , \frac{15}{42} , \frac{17}{42} , \frac{19}{42} , \mathbf{\frac{21}{42}} , \frac{23}{42} , \frac{25}{42} , \frac{27}{42} , \frac{29}{42} , \frac{31}{42} , \frac{33}{42} , \frac{37}{42} , \frac{41}{42} , \frac{45}{42}}$ \\ \hline
    $E_{20}$ & $\scriptscriptstyle{\frac{-5}{66} , \frac{1}{66} , \frac{7}{66} , \frac{13}{66} , \frac{17}{66} , \frac{19}{66} , \frac{23}{66} , \frac{25}{66} , \frac{29}{66} , \frac{31}{66} , \frac{35}{66} , \frac{37}{66} , \frac{41}{66} , \frac{43}{66} , \frac{47}{66} , \frac{49}{66} , \frac{53}{66} , \frac{59}{66} , \frac{65}{66} , \frac{71}{66}}$ \\ \hline
    $J_{4,p}$ & $\scriptscriptstyle{\left\{\frac{12}{12(p+12)}  ,\frac{24}{12(p+12)}  ,\cdots , \frac{12(p+11)}{12(p+12)}\right\}} \textstyle{\cup}$\linebreak $\scriptscriptstyle{\left\{\frac{-(p+12)}{12(p+12)} , \frac{0}{12(p+12)} , \frac{3(p+12)}{12(p+12)} , \frac{4(p+12)}{12(p+12)}  , \frac{5(p+12)}{12(p+12)} ,\frac{6(p+12)}{12(p+12)} , \frac{7(p+12)}{12(p+12)} , \frac{8(p+12)}{12(p+12)}  , \frac{9(p+12)}{12(p+12)} , \frac{12(p+12)}{12(p+12)} , \frac{13(p+12)}{12(p+12)}\right\}}$  \\ \hline
    $E_{24}$ & $\scriptscriptstyle{\frac{-7}{78} , \frac{-1}{78} , \frac{5}{78} , \frac{11}{78} , \frac{17}{78} , \frac{19}{78} , \frac{23}{78} , \frac{25}{78} , \frac{29}{78} , \frac{31}{78} , \frac{35}{78} , \frac{37}{78} , \frac{41}{78} , \frac{43}{78} , \frac{47}{78} , \frac{49}{78} , \frac{53}{78} , \frac{55}{78} , \frac{59}{78} , \frac{61}{78} , \frac{67}{78} , \frac{73}{78} , \frac{79}{78} , \frac{85}{78}}$ \\ \hline
    $E_{25}$ & $\scriptscriptstyle{\frac{-5}{54} , \frac{-1}{54} , \frac{3}{54} , \frac{7}{54} , \frac{11}{54} , \frac{13}{54} , \frac{15}{54} , \frac{17}{54} , \frac{19}{54} , \frac{21}{54} , \frac{23}{54} , \frac{25}{54} , \mathbf{\frac{27}{54}} , \frac{29}{54} , \frac{31}{54} , \frac{33}{54} , \frac{35}{54} , \frac{37}{54} , \frac{39}{54} , \frac{41}{54} , \frac{43}{54} , \frac{47}{54} , \frac{51}{54} , \frac{55}{54} , \frac{59}{54}}$ \\ \hline
    $X_{1,p}$ & $\scriptscriptstyle{\left\{\frac{0}{4(p+4)} , \frac{4}{4(p+4)} ,\cdots , \frac{4(p+4)}{4(p+4)}\right\}}$ $\cup$ $\scriptscriptstyle{\left\{\frac{p+4}{4(p+4)} , \frac{2(p+4)}{4(p+4)} , \frac{2(p+4)}{4(p+4)} , \frac{3(p+4)}{4(p+4)}\right\}}$ \\ \hline
    $Z_{11}$ & $\scriptscriptstyle{\frac{-1}{30} , \frac{5}{30} , \frac{7}{30} , \frac{11}{30} , \frac{13}{30} , \mathbf{\frac{15}{30}} , \frac{17}{30} , \frac{19}{30} , \frac{23}{30} , \frac{25}{30} , \frac{31}{30}}$  \\ \hline
    $Z_{12}$ & $\scriptscriptstyle{\frac{-1}{22} , \frac{3}{22} , \frac{5}{22} , \frac{7}{22} , \frac{9}{22} , \frac{11}{22} , \frac{11}{22} , \frac{13}{22} , \frac{15}{22} , \frac{17}{22} , \frac{19}{22} , \frac{23}{22}}$ \\ \hline
    $Z_{13}$ & $\scriptscriptstyle{\frac{-1}{18} , \frac{2}{18} , \frac{4}{18} , \frac{5}{18} , \frac{7}{18} , \frac{8}{18} , \mathbf{\frac{9}{18}} , \frac{10}{18} , \frac{11}{18} , \frac{13}{18} , \frac{14}{18} , \frac{16}{18} , \frac{19}{18}}$\\ \hline
    $Z_{1,p}$ & $\scriptscriptstyle{\left\{\frac{7}{14(p+7)} , \frac{14}{14(p+7)} , \cdots , \frac{7(2p+13)}{14(p+7)}\right\}}$ $\cup$ $\scriptscriptstyle{\left\{\frac{-(p+7)}{14(p+7)} , \frac{3(p+7)}{14(p+7)}, \frac{5(p+7)}{14(p+7)} , \frac{7(p+7)}{14(p+7)} , \frac{7(p+7)}{14(p+7)} , \frac{9(p+7)}{14(p+7)}, \frac{11(p+7)}{14(p+7)}, \frac{15(p+7)}{14(p+7)}\right\}}$ \\ \hline
    $Z_{17}$ & $\scriptscriptstyle{\frac{-2}{24} , \frac{1}{24} , \frac{4}{24} , \frac{5}{24} , \frac{7}{24} , \frac{8}{24} , \frac{10}{24} , \frac{11}{24} , \mathbf{\frac{12}{24}} , \frac{13}{24} , \frac{14}{24} , \frac{16}{24} , \frac{17}{24} , \frac{19}{24} , \frac{20}{24} , \frac{23}{24} , \frac{26}{24}}$  \\ \hline
    $Z_{18}$ & $\scriptscriptstyle{\frac{-3}{34} , \frac{1}{34} , \frac{5}{34} , \frac{7}{34} , \frac{9}{34} , \frac{11}{34} , \frac{13}{34} , \frac{15}{34} , \frac{17}{34} , \frac{17}{34} , \frac{19}{34} , \frac{21}{34} , \frac{23}{34} , \frac{25}{34} , \frac{27}{34} , \frac{29}{34} , \frac{33}{34} , \frac{37}{34}}$ \\ \hline
    $Z_{19}$ & $\scriptscriptstyle{\frac{-5}{54} , \frac{1}{54} , \frac{7}{54} , \frac{11}{54} , \frac{13}{54} , \frac{17}{54} , \frac{19}{54} , \frac{23}{54} , \frac{25}{54} , \mathbf{\frac{27}{54}} , \frac{29}{54} , \frac{31}{54} , \frac{35}{54} , \frac{37}{54} , \frac{41}{54} , \frac{43}{54} , \frac{47}{54} , \frac{53}{54} , \frac{59}{54}}$  \\ \hline
    $Z_{2,p}$ & $\scriptscriptstyle{\left\{\frac{0}{10(p+10)} , \frac{10}{10(p+10)}, \cdots , \frac{10(p+10)}{10(p+10)}\right\}}\textstyle{\cup}$\linebreak $\scriptscriptstyle{\left\{\frac{-(p+10)}{10(p+10)} , \frac{2(p+10)}{10(p+10)}, \frac{3(p+10)}{10(p+10)}, \frac{4(p+10)}{10(p+10)}, \frac{5(p+10)}{10(p+10)},\frac{5(p+10)}{10(p+10)} , \frac{6(p+10)}{10(p+10)}, \frac{7(p+10)}{10(p+10)}, \frac{8(p+10)}{10(p+10)} , \frac{11(p+10)}{10(p+10)}\right\}}$ \\ \hline
    $Z_{23}$ & $\scriptscriptstyle{\frac{-7}{66} , \frac{-1}{66} , \frac{5}{66} , \frac{11}{66} , \frac{13}{66} , \frac{17}{66} , \frac{19}{66} , \frac{23}{66} , \frac{25}{66} , \frac{29}{66} , \frac{31}{66} , \mathbf{\frac{33}{66}} , \frac{35}{66} , \frac{37}{66} , \frac{41}{66} , \frac{43}{66} , \frac{47}{66} , \frac{49}{66} , \frac{53}{66} , \frac{55}{66} , \frac{61}{66} , \frac{67}{66} , \frac{73}{66}}$  \\ \hline
    $Z_{24}$ & $\scriptscriptstyle{\frac{-5}{46} , \frac{-1}{46} , \frac{3}{46} , \frac{7}{46} , \frac{9}{46} , \frac{11}{46} , \frac{13}{46} , \frac{15}{46} , \frac{17}{46} , \frac{19}{46} , \frac{21}{46} , \frac{23}{46} , \frac{23}{46} , \frac{25}{46} ,\frac{27}{46} , \frac{29}{46} , \frac{31}{46} , \frac{33}{46} , \frac{35}{46} , \frac{37}{46} , \frac{39}{46} , \frac{43}{46} , \frac{47}{46} , \frac{51}{46}}$\\ \hline
    $Z_{25}$ & $\scriptscriptstyle{\frac{-4}{36} , \frac{-1}{36} , \frac{2}{36} , \frac{5}{36} , \frac{7}{36} , \frac{8}{36} , \frac{10}{36} , \frac{11}{36} , \frac{13}{36} , \frac{14}{36} , \frac{16}{36} , \frac{17}{36} , \mathbf{\frac{18}{36}} , \frac{19}{36} , \frac{20}{36} , \frac{22}{36} , \frac{23}{36} , \frac{25}{36} , \frac{26}{36} , \frac{28}{36} , \frac{29}{36} , \frac{31}{36} , \frac{34}{36} , \frac{37}{36} , \frac{40}{36}}$  \\
    \hline
    \end{tabular}
\end{center}
\bigskip

\newpage

\small

\begin{landscape}
\noindent \underline{Adjacency tree to $Z_{3,0}$:}
\bigskip

\def\objectstyle{\scriptstyle}
\def\labelstyle{\scriptstyle}
\xymatrix@= .6pc{
A_1&A_2\ar[l]^-{}&A_3\ar[l]^-{}&A_4\ar[l]^-{}&A_5\ar[l]^-{}&A_6\ar[l]^-{}&A_7\ar[l]^-{}&A_8\ar[l]^-{}&A_9\ar[l]^-{}&A_{10}\ar[l]^-{}&A_{11}\ar[l]^-{}&A_{12}\ar[l]^-{}&&&&&&&&&&&&&\\
&&&D_4\ar[ul]^-{}&D_5\ar[ul]^-{}\ar[l]^-{}&D_6\ar[ul]^-{}\ar[l]^-{}&D_7\ar[ul]^-{}\ar[l]^-{}&D_8\ar[ul]^-{}\ar[l]^-{}&D_9\ar[ul]^-{}\ar[l]^-{}&D_{10}\ar[ul]^-{}\ar[l]^-{}&D_{11}\ar[ul]^-{}\ar[l]^-{}&D_{12}\ar[ul]^-{}\ar[l]^-{}&D_{13}\ar[ul]^-{}\ar[l]^-{}&D_{14}\ar[l]^-{}&&&&&&&&&&&\\
 & & & & &E_6\ar[uul]|\hole^-{}\ar[ul]^-{} &E_7\ar[uul]|\hole^-{}\ar[ul]^-{}\ar[l]^-{} &E_8\ar[uul]|\hole^-{}\ar[ul]^-{}\ar[l]^-{} & & & & & & & & & & & & & & & & & \\
 & & & & & & & & &J_{2,0}\ar[ull]^-{} &J_{2,1}\ar[l]^-{} &J_{2,2}\ar[l]^-{} &J_{2,3}\ar[l]^-{} &J_{2,4}\ar[l]^-{} &J_{2,5}\ar[l]^-{} &J_{2,6}\ar[l]^-{} &J_{2,7}\ar[l]^-{} & & & & & & & & \\
 & & & & & & & & & & &E_{12}\ar[ul]^-{} &E_{13}\ar[ul]^-{}\ar[l]^-{} &E_{14}\ar[ul]^-{}\ar[l]^-{} & & & & & & & & & & & \\
 & & & & & & & & & & & & & & &J_{3,0}\ar[ull]^-{} &J_{3,1}\ar[l]^-{} &J_{3,2}\ar[l]^-{} &J_{3,3}\ar[l]^-{} &J_{3,4}\ar[l]^-{} & & & & & \\
 & & & & & & & & & & & & & & & & &E_{18}\ar[ul]^-{} &E_{19}\ar[ul]^-{}\ar[l]^-{} &E_{20}\ar[ul]^-{}\ar[l]^-{} & & & & & \\
 & & & & & & & & & & & & & & & & & & & & &J_{4,0}\ar[ull]^-{} &J_{4,1}\ar[l]^-{} &J_{4,2}\ar[l]^-{} & \\
 & & & & & & & & & & & & & & & & & & & & & & &E_{24}\ar[ul]^-{} &E_{25}\ar[ul]^-{}\ar[l]^-{} \\
 & & & & & & & &X_{1,0}\ar[uuuuuuull]^-{} &X_{1,1}\ar[l]^-{} &X_{1,2}\ar[l]^-{} &X_{1,3}\ar[l]^-{} &X_{1,4}\ar[l]^-{} &X_{1,5}\ar[l]^-{} &X_{1,6}\ar[l]^-{} &X_{1,7}\ar[l]^-{} &X_{1,8}\ar[l]^-{} &X_{1,9}\ar[l]^-{} & & & & & & & \\
 & & & & & & & & & &Z_{11}\ar[ul]^-{} &Z_{12}\ar[ul]^-{}\ar[l]^-{} &Z_{13}\ar[ul]^-{}\ar[l]^-{} & & & & & & & & & & & & \\
 & & & & & & & & & & & & & &Z_{1,0}\ar[ull]^-{} &Z_{1,1}\ar[l]^-{} &Z_{1,2}\ar[l]^-{} &Z_{1,3}\ar[l]^-{} &Z_{1,4}\ar[l]^-{} &Z_{1,5}\ar[l]^-{} &Z_{1,6}\ar[l]^-{} & & & & \\
 & & & & & & & & & & & & & & & &Z_{17}\ar[ul]^-{} &Z_{18}\ar[ul]^-{}\ar[l]^-{} &Z_{19}\ar[ul]^-{}\ar[l]^-{} & & & & & & \\
 & & & & & & & & & & & & & & & & & & & &Z_{2,0}\ar[ull]^-{} &Z_{2,1}\ar[l]^-{} &Z_{2,2}\ar[l]^-{} &Z_{2,3}\ar[l]^-{} &Z_{2,4}\ar[l]^-{} \\
 & & & & & & & & & & & & & & & & & & & & & &Z_{23}\ar[ul]^-{} &Z_{24}\ar[ul]^-{}\ar[l]^-{} &Z_{25}\ar[ul]^-{}\ar[l]^-{} \\
}
\end{landscape}

\bibliographystyle{unsrt}

\end{document}